\documentclass[twocolumn]{autart}

\usepackage{mathptmx} 

\usepackage{cite}
\usepackage{amsfonts}
\usepackage{algorithmic}
\usepackage{enumerate}
\usepackage{url}
\usepackage{textcomp}

\usepackage{amsmath}
\usepackage{amssymb}
\usepackage{mathtools}
\usepackage{float}
\usepackage{graphicx}
\usepackage{subfigure}
\usepackage{epstopdf}
\usepackage{lineno,hyperref}
\usepackage{xcolor}
\usepackage{multirow}

\newtheorem{definition}{Definition}
\newtheorem{remark}{Remark}
\newtheorem{theorem}{Theorem}

\definecolor{ForestGreen}{rgb}{0.3, 0.7, 0.3}

\modulolinenumbers[5]

\journal{arxiv.org}

\begin{document}

\begin{frontmatter}

\title{\textbf{\textcolor{black}{Event-Triggered Extremum Seeking Control Systems}\thanksref{footnoteinfo}}}

\thanks[footnoteinfo]{\textcolor{black}{Preliminary versions of this paper were partially presented at the 14th IFAC International Workshop on Adaptive and Learning Control Systems \cite{ALCOS:2022}, Casablanca, 2022, and the 22nd World Congress of the International Federation of Automatic Control \cite{IFAC:2023}, Yokohama, 2023.}}


\author[ufrj]{Victor Hugo Pereira Rodrigues},
\author[uerj]{Tiago Roux Oliveira},\linebreak
\author[ufrj]{Liu Hsu},
\author[ucsd]{Mamadou Diagne},
\author[ucsd]{Miroslav Krstic}

\address[ufrj]{Dept. of Electrical Engineering, Federal University of Rio de Janeiro (UFRJ), RJ 21945-970, Brazil.}
\ead{rodrigues.vhp@gmail.com (Victor Hugo Pereira Rodrigues)} 
\ead{liu@coep.ufrj.br (Liu Hsu)}
\address[uerj]{Dept. of Electronics and Telecommunication Engineering, State University of Rio de Janeiro (UERJ), RJ 20550-900, Brazil.}
\ead{tiagoroux@uerj.br (Tiago Roux Oliveira)} 
\address[ucsd]{Dept. of Mechanical and Aerospace Engineering, University of California - San Diego (UCSD), La Jolla, CA 92093-0411, USA.}
\ead{mdiagne@eng.ucsd.edu (Mamadou Diagne)} 
\ead{krstic@ucsd.edu (Miroslav Krstic)}

%
%


\begin{abstract}
This paper proposes an event-triggered control scheme for multivariable extremum seeking of static maps. Both static and dynamic triggering conditions are developed.  Integrating Lyapunov and averaging theories for discontinuous systems, a systematic  design procedure and stability analysis are developed. Both event-based methods enable one to achieve an asymptotic stability result. Ultimately, the resulting closed-loop dynamics demonstrates the advantages of combining both approaches, namely, event-triggered control and extremum seeking. \textcolor{black}{Although
we keep the presentation using the classical event-triggered method, the extension of the results for the periodic event-triggered approach is also indicated.} An illustration of the benefits of the new control method is presented using consistent simulation results, which compare the static and the dynamic triggering approaches.
\end{abstract}
%
%
\end{frontmatter}


\allowdisplaybreaks
\interdisplaylinepenalty=2500 


\section{Introduction}

In spite of the fact that the concept of extremum seeking control (ESC) was introduced 100 years ago by the French Engineer Maurice Leblanc already in 1922 \cite{L:1922}, a rigorous demonstration of stability under feedback control only appeared about 20 years ago \cite{KW:2000}. ESC based on the perturbation method 
adds a periodic excitation-dither signal of small amplitude to the input of the nonlinear map and estimates the gradient by using a suitable demodulation process with the same kind of periodic signals. Extremum seeking is a control strategy that allows the output of a nonlinear map to be held within a vicinity of its extremum. When the parameters of the nonlinear map are available, it is possible to obtain exactly the gradient of the nonlinearity and the control objective can be defined as its stabilization. However, because of parametric uncertainties, the gradient is not always known and this control task is not always straightforward. Despite the several ESC strategies found in the literature, the methods based on perturbations (dither signals) are the oldest and, even remain nowadays quite popular.

After the consolidation of ESC stability results for static and general nonlinear dynamic systems in continuous time \cite{KW:2000}, discrete-time systems \cite{CKAL:2002}, stochastic systems \cite{MK:2009,LK:2010}, multivariable systems \cite{GKN:2012}, and non-cooperative games \cite{FKB:2012}, the theoretical advances of ESC overcome the border of finite-dimension systems to arrive in the world of Partial Differential Equations (PDEs). For such infinite-dimensional systems, boundary feedback controllers are designed to ensure convergence of the ESC in the closed-loop form by compensating the input and/or output under propagation through transport PDEs (delays)  \cite{OKT:2017,REOGK:2019}, 
diffusion PDEs \cite{FKKO:2018,REOGK:2020}, wave PDEs \cite{OK:2021}, Lighthill–Whitham–Richards PDEs \cite{YKOK:2021}, 
parabolic PDEs \cite{REOGK:2021}, and PDE-PDE cascades \cite{OK:2021}.   

In the current technological age of network science, researchers are focusing on decreasing  wiring costs by designing  fast and reliable communication schemes where the plant and controller might not be physically connected, or might even be in different geographical locations. These networked control systems  offer advantages in the financial cost of installation and maintenance \cite{ZHGDDYP:2020}. However, one of their major disadvantages is the resulting high-traffic congestion, which can lead to transmission delays and  packet dropouts, {\it i.e.}, data may be lost while in transit through the network \cite{HNX:2007}. These issues are highly related to the limited resource or available communication channels' bandwidth. To alleviate or mitigate this problem, Event-Triggered Controllers (ETC) can be used.

ETC executes the control task, non-periodically, in response to a triggering condition designed as a function of the plant's state. 
Besides the asymptotic stability properties \cite{T:2007}, this strategy reduces control effort since the control update and data communication only occur when the error between the current state and the equilibrium set exceeds a value that might induce instability \cite{BH:2013}. Pioneering works towards the development of resource-aware control design includes the construction of digital computer design  \cite{s9b}, the event-based PID design  \cite{s9}  and the event-based controller for stochastic systems \cite{s8}. Works dedicated to the extension of event-based control for   networked control systems with a high level of complexity exist  for both linear \cite{s1,HJT:2012,s5} and nonlinear systems \cite{T:2007,APDN:2016}. 
Among others, results on event-based control  deal with the robustness against the effect of  possible  perturbations \cite{s13,s14} or parametric uncertainties \cite{s16}. In  \cite{ZLJ:2021},  ETC is designed to satisfy a cyclic-small-gain condition such that the stabilization of a class of nonlinear time-delay systems is guaranteed. As well, the authors in \cite{CL:2019} proposed  distributed event-triggered leaderless and leader-follower consensus control laws for general linear multi-agent networks. An event-triggered output-feedback  design \cite{APDN:2016} has been employed aiming to stabilize a class of nonlinear systems by combining techniques from event-triggered and time-triggered control. Recently, substantial works have been carried out to conceive event-based approaches for infinite-dimensional systems \cite{s18,s25a,s22,diagne2021event}. 
We emphasize that among existing results \cite{RDEK:2021,rathnayake2022sampled} are focused on infinite-dimensional observer-based event-triggered control for reaction-diffusion PDEs.

This paper introduces a systematic methodology for multivariable event-triggered extremum seeking based on the classical periodic-perturbation method \cite{KW:2000}. \textcolor{black}{Previous conferences versions focusing on the scalar case were presented in \cite{ALCOS:2022} and \cite{IFAC:2023}.} \textcolor{black}{Both topologies for control activation named classical event-triggered and periodic event-triggered mechanisms are explored \cite{s5}.} We consider the design and analysis for static multi-input maps within both static and dynamic event-triggered control frameworks. \textcolor{black}{We prove how this type of extremum seeking approach can deal with closed-loop feedback dynamics that exhibit features typical of 
event-triggered control \cite{T:2007}. ETC---or even periodic ETC \cite{s5}---differs from standard periodic sampled-data control \cite{KNTM:2013,HNW:2023,ZFO:2023}, as in ETC the event times (which result from the triggering condition and the system's state evolution) are in general only a (specific) subset of the sampling times and can be aperiodic.} The stability analysis is carried out using time-scaling technique, averaging theory for discontinuous systems and Lyapunov's direct method to characterize the properties of the closed-loop system. 
\textcolor{black}{
Note that the non-periodic nature of the ETC approach is not an obstacle to apply the averaging result for discontinuous systems by Plotnikov \cite{P:1979}---see Appendix~\ref{appendix_plotnikov}---since the perturbation-probing signals of the ESC are still periodic functions. Conversely, the discontinuities are not in the periodic perturbations, but in the states. In addition, the event-triggered extremum seeking invention is totally different from the previous related literature which considers, for instance, dead-zones \cite{deadzone_2016} or time-delays \cite{OKT:2017}. Although artificial delays are able to retard the system response and dead zones are a common way to automatically turn off a controller when the system's
error is below a prescribed threshold, none of them can compete with ETC capabilities in order to preserve control performance under limited  bandwidth of the actuation-sensing paths and still guarantee stability for the closed-loop dynamics. Particularly, as shown in \cite{OKT:2017}, arbitrary time-delays may destabilize ESC systems.}

This paper is organized as follows:  Section~\ref{sec:prblFrm} presents the formulation of the control problem. For both static and dynamic, the triggering conditions are presented in Section~\ref{sec:cls}. Stability and convergence results of multivariable \textit{static} event-triggered extremum seeking are presented in Section~\ref{ETESNC_unknownH*}. Stability and convergence results of multivariable \textit{dynamic} event-triggered extremum seeking are presented in Section~\ref{DETESNC_unknownH*}. \textcolor{black}{Section~\ref{PETC} brings the alternative implementation using the  periodic event-triggered extremum seeking approach.} Simulation results are discussed in Section~\ref{sec:sim} and concluding remarks in Section~\ref{sec:concl}.

\emph{Notation:} 
Throughout the manuscript, the 2-norm (Euclidean) of vectors and induced norm of matrices are denoted by double bars $\|\cdot\|$ while absolute value of scalar variables are denoted by single bars $|\cdot|$. The terms $\lambda_{\min}(\cdot)$ and $\lambda_{\max}(\cdot)$ denote  the minimum and maximum eigenvalues of a matrix, respectively. Consider $\varepsilon \in \lbrack -\varepsilon_{0}\,, \varepsilon_{0} \rbrack \subset \mathbb{R}$ and the mappings $\delta_{1}(\varepsilon)$ and $\delta_{2}(\varepsilon)$, where $\delta_{1}: \lbrack -\varepsilon_{0}\,, \varepsilon_{0} \rbrack \to \mathbb{R}$ and $\delta_{2}: \lbrack -\varepsilon_{0}\,, \varepsilon_{0} \rbrack \to \mathbb{R}$. One states that $\delta_{1}(\varepsilon)$ has magnitude order of $\delta_{2}(\varepsilon)$, {\it i.e.}, $\delta_{1}(\varepsilon) = \mathcal{O}(\delta_{2}(\varepsilon))$, if there exist positive constants $k$ and $c$ such that $|\delta_{1}(\varepsilon)| \leq k |\delta_{2}(\varepsilon)|$, for all $|\varepsilon|<c$.

\section{Problem Formulation} \label{sec:prblFrm}

\begin{figure*}[h!]
\begin{center}
\includegraphics[width=17.0cm]{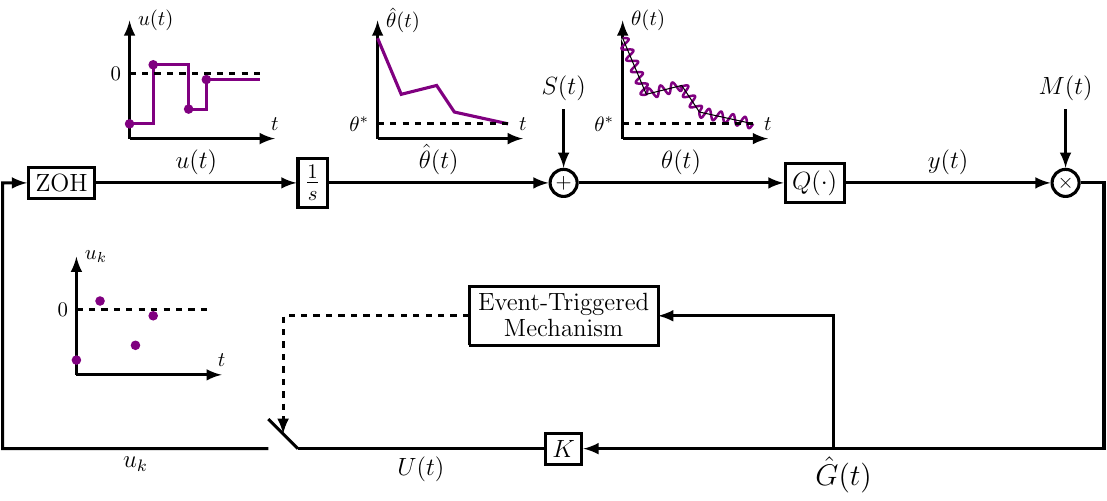}
\end{center}
\caption{\textcolor{black}{Event-triggered based on extremum seeking scheme. The key distinction is the piecewise linearity of $\hat{\theta}(t)$. This comes from the fact that ESC dynamics is essentially an integrator and the fact that the integrator’s input $u(t)$ is being made piecewise-constant by ETC.  Hence, the piecewise linearity of $\hat{\theta}(t)$, with a superimposition of the sinusoidal perturbation $S(t)$, is the one feature of ETC-ESC that is worth highlighting in the input $\theta(t)$ of the nonlinear map $y(t)=Q(\theta(t))$ to be optimized.}}
\label{fig:blockDiagram_4}
\end{figure*}

We define the following  nonlinear static map
\begin{align}
Q(\theta(t)) = Q^{\ast}+\frac{1}{2}(\theta(t)-\theta^{\ast})^{T}H^{\ast}(\theta(t)-\theta^{\ast})\,, \label{eq:Q_1_event}
\end{align}
where $H^{\ast}=H^{\ast T} \in \mathbb{R}^{n \times n}$ is the Hessian matrix, $\theta^{\ast} \in \mathbb{R}^{n}$ is the unknown optimizer, $\theta(t)\in \mathbb{R}^{n}$ is the input map, designed as the real-time estimate $\hat{\theta}(t) \in \mathbb{R}^{n}$ of $\theta^{\ast}$ additively perturbed by the vector $S(t)$ of sinusoids, {\it i.e.},
\begin{align}
\theta(t)=\hat{\theta}(t)+S(t)\,. \label{eq:theta_event}
\end{align}
The output of the nonlinear map (\ref{eq:Q_1_event}) can be written as $y(t)=Q(\theta(t)).$ Fig.~\ref{fig:blockDiagram_4} shows the closed-loop structure of the event-triggered-based extremum seeking control system to be designed. 
\textcolor{black}{Equation (\ref{eq:Q_1_event}) states that, if maps are at least locally quadratic, our methodology can be applied, and provide guarantees, in some neighborhood of the extremum. For maps that are not locally quadratic, and may not yield exponential stability of the
average system, an approach by \cite{TNM:2006} might lead to \textit{asymptotic} practical stability, not exponential, but the form of averaging theory that they use is not available for systems on Banach spaces or closed-loop dynamics with discontinuous right-hand sides.}

The following assumptions are considered throughout the paper.
\begin{assum}\label{assumption1}
The unique optimizer vector $\theta^{\ast} \in \mathbb{R}^{n}$, the Hessian matrix $H^{\ast}$ and the scalar $Q^{\ast} \in \mathbb{R}$ are unknown parameters of the nonlinear map (\ref{eq:Q_1_event}). Moreover, $H^{\ast}$ is symmetric and has known definite sign, thus, being full rank.
\end{assum}

\begin{assum}\label{assumption2}
The matrix product $H^{\ast}K$ is Hurwitz such that for any given $Q=Q^{T}>0$ there exists a $P=P^{T}>0$ that satisfies the Lyapunov equation
\begin{align}
K^{T}H^{\ast T}P+PH^{\ast}K=-Q \,. \label{eq:LyapEq}
\end{align}
Furthermore, the sum of the induced norms of the matrices $K^{T}H^{\ast T}P$ and $PH^{\ast}K$ is upper bounded by a known positive constant $\beta$,
\begin{align}
\|K^{T}H^{\ast T}P\|+\|PH^{\ast}K\|\leq \beta \,, \label{eq:b}
\end{align}
and
the least eigenvalue of the matrix $Q$ is lower bounded by a known positive constant $\alpha$,
\begin{align}
\alpha \leq \lambda_{\min}(Q)  \,. \label{eq:a}
\end{align}
\end{assum}
\vspace{-0.25cm}
\textcolor{black}{
Assumptions~\ref{assumption1} and \ref{assumption2} are not asking for a full-knowledge of the map parameters, thus respecting the spirit of the ESC approach, that is supposed to be used for unknown maps, with unknown parameters $\theta^*$, $Q^*$ and $H^*$. In particular, the symmetry condition imposed to $H^*$ simplifies the problem and will be necessary in the stability analysis.}

\subsection{Continuous-Time Extremum Seeking}

Let us define the estimation error 
\begin{align}
\tilde{\theta}(t)=\hat{\theta}(t)-\theta^{\ast}\,, \label{eq:thetaTilde_event}
\end{align}
and the Gradient estimate by the demodulation
\begin{align}
\hat{G}(t)=M(t)y(t)\,, \label{eq:hatG_event}
\end{align}
with dither vectors (see \cite{GKN:2012,K:2014})
\begin{align}
S(t)&= \left[a_{1}\sin\left(\omega_1 t\right),\ldots\,,a_{n}\sin\left(\omega_n t\right)\right]^{T} \,, \label{eq:S_event} \\
M(t)&=\left[\frac{2}{a_{1}}\sin\left(\omega_1 t\right),\ldots\,,\frac{2}{a_{n}}\sin\left(\omega_n t\right)\right]^{T} \,, \label{eq:M_event}
\end{align} 
of nonzero amplitudes $a_{i}$. Moreover, the probing frequencies $\omega_{i}$'s can be selected as
\begin{align}
\omega_{i}=\omega_{i}'\omega \,, \quad i \in \left\{1,\ldots\,,n\right\}\,, \label{eq:omegai_event}
\end{align}
where $\omega$ is a positive constant and $\omega_{i}'$ is a rational number.

\begin{assum}\label{assumption3}
The probing frequencies satisfy
\begin{align}
\omega'_{i} 	\notin \left\{\omega'_{j}\,,~\frac{1}{2}(\omega'_{j}+\omega'_{k})\,,~\omega'_{j}+2\omega'_{k}\,,~\omega'_{k}\pm\omega'_{l}\right\}\,, \label{eq:omega_iNotIn}
\end{align}
for all $i$, $j$, $k$ and $l$.
\end{assum} 

From (\ref{eq:theta_event}) and (\ref{eq:thetaTilde_event}), one has
\begin{align}
\theta(t)=\tilde{\theta}(t)+S(t)+\theta^{\ast}\,, \label{eq:theta_2_event}
\end{align}
and, therefore, by plugging (\ref{eq:theta_2_event}) into (\ref{eq:Q_1_event}), $y(t)$ can be written as
\begin{align}
y(t)&=Q^{\ast}+\frac{1}{2}(\tilde{\theta}(t)+S(t))^{T} H^{\ast}(\tilde{\theta}(t)+S(t)) \nonumber \\
&=Q^{\ast}+\frac{1}{2}\tilde{\theta}^{T}(t)H^{\ast}\tilde{\theta}(t)+S^{T}(t)H^{\ast}\tilde{\theta}(t)\nonumber \\
&\quad+\frac{1}{2}S^{T}(t)H^{\ast}S(t)\,.
 \label{eq:y_1_event}
\end{align}
 Thus, from (\ref{eq:hatG_event}), (\ref{eq:y_1_event}) and following \cite{GKN:2012}, the gradient estimate  in the average sense, is given by
\begin{align}
\hat{G}(t)&= M(t)Q^{\ast}+\frac{1}{2}M(t)\tilde{\theta}^{T}(t)H^{\ast}\tilde{\theta}(t)\nonumber \\
&\quad+M(t)S^{T}(t)H^{\ast}\tilde{\theta}(t)+\frac{1}{2}M(t)S^{T}(t)H^{\ast}S(t)\,. \label{eq:hatG_2_event}
\end{align} 
Now, defining 
\begin{align}
H\left(t\right)&:=M(t)S^{T}(t)H^{\ast} \nonumber \\
&=H^{\ast}+\Delta(t) H^{\ast}\,, \label{eq:hatH_event}
\end{align}
where the elements of the $\Delta (t) \in \mathbb{R}^{n\times n}$ are given by
\begin{align}
\Delta_{ii} (t) &\mathbb{=} -\cos(2\omega_{i}t)\,, \label{eq:Deltaii_event} \\
\Delta_{ij} (t) &\mathbb{=} 2\frac{a_{j}}{a_{i}}\sin(\omega_{i}t)\sin(\omega_{j}t) \nonumber \\
&\mathbb{=}\frac{a_{j}}{a_{i}}\cos((\omega_{i}-\omega_{j})t)\!-\!\frac{a_{j}}{a_{i}}\cos((\omega_{i}+\omega_{j})t)\,,  \label{eq:Deltaij_event}
\end{align}
for all $i\!\neq\! j$, and using \eqref{eq:hatH_event}, one can rewrite (\ref{eq:hatG_2_event}) as follows
\begin{align}
\hat{G}(t)&=H(t)\tilde{\theta}(t)+M(t)Q^{\ast}+\frac{1}{2}H(t)S(t)+\vartheta(t) \,, \label{eq:hatG_3_event} \\
\vartheta(t)&:=\frac{1}{2}M(t)\tilde{\theta}^{T}(t)H^{\ast}\tilde{\theta}(t)\,. \label{eq:vartheta_event}
\end{align} 
The term $\vartheta(t)$ given above is quadratic in $\tilde{\theta}(t)$ and, therefore, may be neglected in a local analysis \cite{K:2014}. 
Thus, hereafter the gradient estimate is given by
\begin{align}
\hat{G}(t)&=H(t)\tilde{\theta}(t)+M(t)Q^{\ast}+\frac{1}{2}H(t)S(t) \,. \label{eq:hatG_4_event}
\end{align}

On the other hand, from the time-derivative of (\ref{eq:thetaTilde_event}) and the ESC scheme depicted in Fig.~\ref{fig:blockDiagram_4}, the dynamics that governs $\hat{\theta}(t)$, as well as $\tilde{\theta}(t)$, is given by
\begin{align}
\dot{\tilde{\theta}}(t)&=\dot{\hat{\theta}}(t)=u(t) \label{eq:dotTildeTheta_1_event}\,, 
\end{align}
where $u$ is an ESC law to be designed.

By taking the time-derivative of (\ref{eq:hatG_4_event}), with the help of  (\ref{eq:hatH_event}) and (\ref{eq:dotTildeTheta_1_event}), one gets the following differential equation
\begin{align}
\dot{\hat{G}}(t)&=f(t,\tilde{\theta}(t),u(t)) \nonumber \\
&=H(t)u(t)+w(t,\tilde{\theta}(t)) \, \label{eq:dotHatG_event},
\end{align}
where
\begin{align}
w(t,\tilde{\theta}(t))&=\dot{\Delta }(t) H^{\ast}\tilde{\theta}(t)+\dot{M}(t)Q^{\ast}+\frac{1}{2}\dot{\Delta }(t)H^{\ast}S(t) \nonumber \\
&\quad+\frac{1}{2}H^{\ast}\dot{S}(t)+\frac{1}{2}\Delta (t)H^{\ast}\dot{S}(t) \label{eq:W_event} \,.
\end{align}
For all $t\geq 0$, the  continuous-time feedback law
\begin{align}
u(t)=K\hat{G}(t) \,, \quad \forall t\geq 0 \label{eq:U_continuous}
\end{align}
\textcolor{black}{is a stabilizing controller for the average version of (\ref{eq:dotHatG_event}), since the gain $K$ is chosen such that $H^{\ast}K$ is Hurwitz. Recalling that $w(t,\tilde{\theta}(t))$ in (\ref{eq:W_event}) has null average, consequently, one can write:
\begin{align}
\dot{\hat{G}}_{\rm{av}}(t)=H^{\ast}K \hat{G}_{\rm{av}}(t)\,, \quad \forall t\geq 0\,. \label{eq:dotHatG_event_average_golden}
\end{align}
}
Our goal is to design a stabilizing controller for the closed-loop system \eqref{eq:dotHatG_event} and \eqref{eq:W_event} in a sampled and hold fashion \cite{APDNH:2018} by  emulating the continuous-time control law (\ref{eq:U_continuous}). Here, the control law is only updated for a given  sequence of time instants $(t_{k})_{k\in\mathbb{N}}$ defined by an event-generator that is constructed to preserve stability and robustness features. More precisely, the execution of the control task is orchestrated by a monitoring mechanism that invokes control updates when the difference between the current value of the output and its previously computed value  at time $t_k$ becomes too large with respect to a constructed triggering condition that needs to be satisfied \cite{HJT:2012}. Note that in a conventional sampled-data implementation, the execution times are distributed equidistantly in time, meaning that $t_{k+1}=t_{k}+ h$, where $h>0$ is a known constant, for all $k\in \mathbb{N}$, while in event-triggered scheme aperiodic sampling may occur.  

\subsection{Emulation of the Continuous-Time Extremum Seeking Design}

Defining the control input for all $t\in[t_{k},t_{k+1}), k\in\mathbb{N}$ as
\begin{align}
u_k=K\hat{G}(t_{k}) \,, \label{eq:U_event}
\end{align}
we introduce the error vector, that is to say  the deviation of the output signal as 
\begin{align}
e(t):=\hat{G}(t_{k})-\hat{G}(t) \,, \quad \forall t \in \lbrack t_{k}\,, t_{k+1}) \,, \quad k\in \mathbb{N} \,. \label{eq:e_event}
\end{align}

Now, using the event-triggered control law (\ref{eq:U_event}), adding and subtracting the term  $H(t)K\hat{G}(t)$ and $K\hat{G}(t)$ to (\ref{eq:dotHatG_event}) and  (\ref{eq:dotTildeTheta_1_event}), respectively, one arrives at the Input-to-State Stable (ISS) representation of \eqref{eq:dotTildeTheta_1_event} and \eqref{eq:dotHatG_event} with respect to the error vector  (\ref{eq:e_event}) and the time-varying disturbance $w(t,\tilde{\theta}(t))$. The resulting dynamics are given below, \textcolor{black}{$\forall t \in \lbrack t_{k}\,, t_{k+1}) \,,  ~k\in \mathbb{N}$}:
\begin{align}
\dot{\hat{G}}(t)&\mathbb{=} H(t)K\hat{G}(t)+H(t)Ke(t)+w(t,\tilde{\theta}(t))\,, \label{eq:dotHatGav_event_3} \\
\dot{\tilde{\theta}}(t)&\mathbb{=} KH(t)\tilde{\theta}(t)+Ke(t)+KM(t)Q^{\ast}\!+\!\frac{1}{2}KH(t)S(t)\,. \label{eq:dotTildeTheta_2_event}
\end{align}
In the subsequent developments, the static and dynamic triggering mechanisms, presented in Definitions~\ref{def:staticEvent} and \ref{def:dynamicEvent}, respectively.

\subsection{Event-Triggered Control Mechanism}

Definitions~\ref{def:staticEvent} and \ref{def:dynamicEvent} show how the nonlinear mapping $\Xi: \mathbb{R}^n \times \mathbb{R}^n \mapsto \mathbb{R}$ , given by 
\begin{align}
\Xi(\hat{G},e)&=\sigma \alpha \|\hat{G}(t)\|^2-\beta \|e(t)\|\|\hat{G}(t)\|	\,, \label{eq:Xi_event_1}
\end{align}
where, $\sigma \in (0,1)$ is a given parameter, can be employed in the design of the static and dynamic execution mechanisms. The  mapping $\Xi(\hat{G},e)$ is designed to appropriately re-compute  the control law (\ref{eq:U_event}) and update the ZOH actuator depicted in Fig.~\ref{fig:blockDiagram_4} such that the asymptotic stability of the closed-loop system is achieved \cite{HJT:2012}.

\begin{definition}[\small{Multivariable Static Triggering Condition}] \label{def:staticEvent}
Let $\Xi(\hat{G},e)$ in (\ref{eq:Xi_event_1}) be the nonlinear mapping  and $K$  the control gain in (\ref{eq:U_event}). The event-triggered controller with static-triggering condition consists of two components:
\begin{enumerate}
	\item A set of increasing sequence of time $I=\{t_{0}\,, t_{1}\,, t_{2}\,,\ldots\}$ with $t_{0}=0$ generated under the following rules:
		\begin{itemize}
			\item If $\left\{t \in\mathbb{R}^{+}: t>t_{k} ~ \wedge ~ \Xi (\hat{G},e) < 0\right\} = \emptyset$, then the set of the times of the events is $I=\{t_{0}\,, t_{1}\,, \ldots, t_{k}\}$.
			\item If $\left\{t \in\mathbb{R}^{+}: t>t_{k} ~ \wedge ~ \Xi (\hat{G},e) < 0\right\} \neq \emptyset$, next event time is given by
				\begin{align}
					t_{k+1}&=\inf\left\{t \in\mathbb{R}^{+}: t>t_{k} ~ \wedge ~ \Xi (\hat{G},e) <0 \right\}\,, \label{eq:tk+1_event}
				\end{align}
				consisting of the static event-trigger mechanism.
		\end{itemize}
	\item A feedback control action updated at the  triggering instants \eqref{eq:U_event}.
\end{enumerate}  
\end{definition}

Although the aperiodicity of the control update is guaranteed by the static event generation mechanism (\ref{eq:tk+1_event}), it is often convenient to use its filtered version  to increase the inter-execution times. In this case, inspired by \cite{G:2014}, we also construct a dynamic event-triggering mechanism described in the following definition.

\begin{definition}[\small{Multivariable Dynamic Triggering Condition}] \label{def:dynamicEvent}
Let $\Xi(\hat{G},e)$ in (\ref{eq:Xi_event_1}) be the nonlinear mapping  and $K$  the control gain in (\ref{eq:U_event}), $\gamma>0$ a positive constant and $\upsilon(t)$ the solution of the dynamics 
\begin{align}
\dot{\upsilon}(t)&=-\mu \upsilon(t) +\Xi (\hat{G},e) \,, \quad \mu>0\,, \quad \upsilon(0)\geq 0\,. \label{eq:eta}
\end{align}
The event-triggered controller with dynamic triggering condition consists of two components:
\begin{enumerate}
	\item A set of increasing sequence of time $I=\{t_{0}\,, t_{1}\,, t_{2}\,,\ldots\}$ with $t_{0}=0$ generated under the following rules:
		\begin{itemize}
			\item If $\left\{t \in\mathbb{R}^{+}: t>t_{k} ~ \wedge ~ \upsilon(t)+\gamma\Xi (\hat{G},e) < 0\right\} = \emptyset$, then the set of the times of the events is $I=\{t_{0}\,, t_{1}\,, \ldots, t_{k}\}$.
			\item If $\left\{t \in\mathbb{R}^{+}: t>t_{k} ~ \wedge ~ \upsilon(t)+\gamma\Xi (\hat{G},e) < 0\right\} \neq \emptyset$, next event time is given by
				\begin{align}
					\!\!\!\!\!\!\!\!\! t_{k+1}&=\inf\left\{t \in\mathbb{R}^{+}: t>t_{k} ~ \wedge ~ \upsilon(t)+\gamma\Xi (\hat{G},e) <0 \right\}\!, \label{eq:tk+1_event_dynamic}
				\end{align}
				which is the dynamic event-trigger mechanism.
		\end{itemize}
	\item A feedback control action updated at the triggering instants \eqref{eq:U_event}.
\end{enumerate}  
For all $t\in(t_{k},t_{k+1})$, $\upsilon(t_{0})=\upsilon(0)\geq 0$ and $\upsilon(t_{k}^{-})=\upsilon(t_{k})=\upsilon(t_{k}^{+})$.
\end{definition}

\section{Closed-loop system and Triggering Conditions} \label{sec:cls}

\subsection{Time-Scaling System}

Now, we  introduce a suitable time scale to carry out  the stability analysis of the closed-loop system. From (\ref{eq:omegai_event}), one can notice that the dither frequencies  (\ref{eq:S_event}) and (\ref{eq:M_event}), as well as their combinations (\ref{eq:Deltaii_event}) and (\ref{eq:Deltaij_event}), are both rational. Furthermore, there exists a time period $T$ such that 
\begin{align}
T&= 2\pi \times \text{LCM}\left\{\frac{1}{\omega_{i}}\right\}\,, \quad \forall i \left\{1\,,2\,,\ldots\,,n\right\}\,, \label{eq:T}
\end{align}
 where LCM denotes the least common multiple such that it is possible to define the time-scale for the dynamics (\ref{eq:dotHatGav_event_3})--(\ref{eq:dotTildeTheta_2_event}) with the transformation $\bar{t}=\omega t$, where 
\begin{align}
\omega&:=\frac{2\pi}{T}\,. \label{eq:omega_event_1}
\end{align}
 Hence, the system (\ref{eq:dotHatGav_event_3}), (\ref{eq:dotTildeTheta_2_event}) and (\ref{eq:eta}) can be rewritten as, \textcolor{black}{$\forall t \in \lbrack t_{k}\,, t_{k+1}) \,,  ~k\in \mathbb{N}$}, 
\begin{align}
\frac{d\hat{G}(\bar{t})}{d\bar{t}}&=\frac{1}{\omega}\hat{\mathcal{G}}\left(\bar{t},\hat{G},\tilde{\theta},\upsilon,\dfrac{1}{\omega}\right)\,, \label{eq:dotHatGav_event_4} \\
\frac{d\tilde{\theta}(\bar{t})}{d\bar{t}}&= \frac{1}{\omega}\tilde{\Theta}\left(\bar{t},\hat{G},\tilde{\theta},\upsilon,\dfrac{1}{\omega}\right)\,, \label{eq:dotTildeTheta_3_event} \\
\frac{d\upsilon(\bar{t})}{d\bar{t}}&= \frac{1}{\omega}\tilde{\Upsilon}\left(\bar{t},\hat{G},\tilde{\theta},\upsilon,\dfrac{1}{\omega}\right)\,, \label{eq:dotUpsilon_}
\end{align}
where
\begin{align}
\hat{\mathcal{G}}\left(\bar{t},\hat{G},\tilde{\theta},\upsilon,\dfrac{1}{\omega}\right)
&=H(\bar{t})K\hat{G}(\bar{t})+H(\bar{t})Ke(\bar{t}) \nonumber \\
&\quad +w(\bar{t},\tilde{\theta}(\bar{t}))\,, \label{eq:hatMathcalG} \\
\tilde{\Theta}\left(\bar{t},\hat{G},\tilde{\theta},\upsilon,\dfrac{1}{\omega}\right)&=KH(\bar{t})\tilde{\theta}(\bar{t})+Ke(\bar{t}) \nonumber \\
&\quad+KM(\bar{t})Q^{\ast}+\frac{1}{2}KH(\bar{t})S(\bar{t})\,, \label{eq:tildeMathcalTheta} \\
\Upsilon\left(\bar{t},\hat{G},\tilde{\theta},\upsilon,\dfrac{1}{\omega}\right)&=-\mu \textcolor{black}{\upsilon(\bar{t})}+\Xi (\hat{G},e) \,.\label{eq:dotUpsilon_3_event}
\end{align}
From the above dynamics, an appropriate averaging system in the new time-scale $\bar{t}$ can be introduced. \textcolor{black}{The right-hand side of (\ref{eq:dotHatGav_event_4})--(\ref{eq:dotUpsilon_}) is
discontinuous, and the discontinuity is not periodic in the states (since the
triggering events are not periodic), but the system is still periodic in time. It is then clear why the averaging results by Plotnikov \cite{P:1979} are applicable to this scheme---see Appendix~\ref{appendix_plotnikov}.}

\subsection{Average System}

Defining the augmented state as follows
\begin{align}
X^{T}(\bar{t}):=\begin{bmatrix} \hat{G}^{T}(\bar{t})\,, \tilde{\theta}^{T}(\bar{t}) \,, \upsilon(\bar{t})
\end{bmatrix}\,,
\end{align}
the system (\ref{eq:dotHatGav_event_4})--(\ref{eq:dotUpsilon_3_event}) reduces to
\begin{align}
\dfrac{dX(\bar{t})}{d\bar{t}}&=\dfrac{1}{\omega}\mathcal{F}\left(\bar{t},X,\dfrac{1}{\omega}\right)\,, \label{eq:dotX_event}
\end{align}
where $\mathcal{F}^{T}=\begin{bmatrix} \hat{\mathcal{G}}^{T}\,, \tilde{\Theta}^{T}\,, \Upsilon \end{bmatrix}$. Note that (\ref{eq:dotX_event}) is characterized by a small parameter $1/\omega$ as well as a $T$-periodic function $\mathcal{F}\left(\bar{t},X,\dfrac{1}{\omega}\right)$ in $\bar{t}$ and, thereby, the averaging method can be performed on  $\mathcal{F}\left(\bar{t},X,\dfrac{1}{\omega}\right)$ at $\displaystyle \lim_{\omega\to \infty}\dfrac{1}{\omega}=0$, as shown in references \cite{K:2002,P:1979}. The averaging method allows for determining in what sense the behavior of a constructed average autonomous system approximates the behavior of the non-autonomous system (\ref{eq:dotX_event}). By employing the averaging technique to (\ref{eq:dotX_event}), we derive the following average system
\begin{align}
\dfrac{dX_{\text{av}}(\bar{t})}{d\bar{t}}&=\dfrac{1}{\omega}\mathcal{F}_{\text{av}}\left(X_{\text{av}}\right) \,, \label{eq:dotXav_event_1} \\
\mathcal{F}_{\text{av}}\left(X_{\text{av}}\right)&=\dfrac{1}{T}\int_{0}^{T}\mathcal{F}\left(\delta,X_{\text{av}},0\right)d\delta
\,,  \label{eq:mathcalFav_event}
\end{align}
where the averaging terms are given below
\begin{align}
S_{\text{av}}(\bar{t})&\mathbb{=} \frac{1}{T}\int_{0}^{T}S(\delta)d\delta=0 \,,  \dot{S}_{\text{av}}(\bar{t})\mathbb{=} \frac{1}{T}\int_{0}^{T}\dot{S}(\delta)d\delta=0\,, \label{eq:Sav_event} \\
M_{\text{av}}(\bar{t})&\mathbb{=} \frac{1}{T}\int_{0}^{T}M(\delta)d\delta=0 \,,  \dot{M}_{\text{av}}(\bar{t})\mathbb{=} \frac{1}{T}\int_{0}^{T}\dot{M}(\delta)d\delta=0\,, \label{eq:Mav_event} \\
\Delta_{\text{av}}(\bar{t})&\mathbb{=} \frac{1}{T}\int_{0}^{T}\Delta(\delta)d\delta\mathbb{=}0\,,  \dot{\Delta}_{\text{av}}(\bar{t})\mathbb{=} \frac{1}{T}\int_{0}^{T}\dot{\Delta}(\delta)d\delta=0\,, \label{eq:Deltaav_event}
\end{align}
and, consequently, 
\begin{align}
H_{\text{av}}(\bar{t})&\mathbb{=}\frac{1}{T}\int_{0}^{T}H(\delta)d\delta=H^{\ast} \,,
\dot{H}_{\text{av}}(\bar{t})\mathbb{=}\frac{1}{T}\int_{0}^{T}\dot{H}(\delta)d\delta
=0 \,. \label{eq:dotHav_event}
\end{align}
\textcolor{black}{Therefore, we ``freeze'' the average states of $\hat{G}(\bar{t})$, $e(\bar{t})$ and $\tilde{\theta}(\bar{t})$} in (\ref{eq:dotHatGav_event_3}) and (\ref{eq:W_event}), and by using the averaging values (\ref{eq:Sav_event})--(\ref{eq:dotHav_event}), one gets for all $\bar t\in [\bar t_k, \bar t_{k+1})$
\begin{align}
\frac{d\hat{G}_{\text{av}}(\bar{t})}{d\bar{t}}&=\frac{1}{\omega}H^{\ast}K\hat{G}_{\text{av}}(\bar{t})+\frac{1}{\omega}H^{\ast}Ke_{\text{av}}(\bar{t})\,, \label{eq:dotHatGav_event_1} \\
e_{\text{av}}(\bar{t})&=\hat{G}_{\text{av}}(\bar{t}_{k})-\hat{G}_{\text{av}}(\bar{t})\,, \label{eq:Eav_event_1} 
\end{align}
since the average value of $w(t,\tilde{\theta}(t))$ in (\ref{eq:W_event}) is
\begin{align}
w_{\text{av}}(\bar{t},\tilde{\theta}_{\text{av}}(\bar{t}))&=\dot{\Delta }_{\text{av}}(\bar{t}) H^{\ast}\tilde{\theta}_{\text{av}}(\bar{t})+\dot{M}_{\text{av}}(\bar{t})Q^{\ast} \nonumber \\
&\quad+\frac{1}{2}\dot{\Delta }_{\text{av}}(\bar{t})H^{\ast}S_{\text{av}}(\bar{t}) \nonumber \\
&\quad+\frac{1}{2}H^{\ast}\dot{S}_{\text{av}}(\bar{t})+\frac{1}{2}\Delta_{\text{av}} (\bar{t})H^{\ast}\dot{S}_{\text{av}}(\bar{t}) \label{eq:Wav_event} \nonumber \\
&=0\,.
\end{align}
Hence, from (\ref{eq:dotHatGav_event_1}) it is easy to verify the ISS relationship of $\hat{G}_{\text{av}}(\bar{t})$ with respect to the averaged measurement error  $e_{\text{av}}(\bar{t})$ in (\ref{eq:Eav_event_1}). 

Moreover, from (\ref{eq:hatG_4_event}), one has
\begin{align}
\hat{G}_{\text{av}}(\bar{t})= H^{\ast}\tilde{\theta}_{\text{av}}(\bar{t})\,, \label{eq:hatGav_event_1}
\end{align}
and, consequently,
\begin{align}
\tilde{\theta}_{\text{av}}(\bar{t})= H^{\ast-1}\hat{G}_{\text{av}}(\bar{t})\,. \label{eq:tildeThetaAv_event_1}
\end{align}
Taking the time-derivative of (\ref{eq:tildeThetaAv_event_1}), we get
\begin{align}
\frac{d\tilde{\theta}_{\text{av}}(\bar{t})}{d\bar{t}}=\frac{1}{\omega}KH^{\ast}\tilde{\theta}_{\text{av}}(\bar{t})+\frac{1}{\omega}Ke_{\text{av}}(\bar{t})\,. \label{eq:dotTildeThetaAv_event_1}
\end{align}

Therefore, the following average event-triggered detection laws can be introduced for the average system.

Defining the average version of $\Xi(\hat{G},e)$, {\it i.e.} as
\begin{align}
\Xi(\hat{G}_{\rm{av}},e_{\rm{av}})&=\sigma \alpha \|\hat{G}_{\rm{av}}(\bar{t})\|^2-\beta \|e_{\rm{av}}(\bar{t})\|\|\hat{G}_{\rm{av}}(\bar{t})\|	\,, \label{eq:Xi_event_2}
\end{align}
we construct the average event-triggered mechanisms.

\begin{definition}[\small{Average Static Triggering Condition}] 
\label{def:averageStaticEvent} Let $\Xi(\hat{G}_{\rm{av}},e_{\rm{av}})$ in (\ref{eq:Xi_event_2}) be the nonlinear mapping  and $K$  the control gain in (\ref{eq:U_event}). The event-triggered controller with average static-triggering condition in the new time-scale consists of two components:
\begin{enumerate}
	\item A set of increasing sequence of time $I=\{\bar{t}_{0}\,, \bar{t}_{1}\,, \bar{t}_{2}\,,\ldots\}$ with $\bar{t}_{0}=0$ generated under the following rule:
		\begin{itemize}
			\item If $\left\{\bar{t} \in\mathbb{R}^{+}: \bar{t}>\bar{t}_{k} ~ \wedge ~ \Xi (\hat{G}_{\rm{av}},e_{\rm{av}}) < 0\right\} = \emptyset$, then the set of the times of the events is $I=\{\bar{t}_{0}\,, \bar{t}_{1}\,, \ldots, \bar{t}_{k}\}$.
			\item If $\left\{\bar{t} \in\mathbb{R}^{+}: \bar{t}>\bar{t}_{k} ~ \wedge ~ \Xi (\hat{G}_{\rm{av}},e_{\rm{av}}) < 0\right\} \neq \emptyset$, next event time is given by
				\begin{align}
					\bar{t}_{k+1}&=\inf\left\{\bar{t} \in\mathbb{R}^{+}: \bar{t}>\bar{t}_{k} ~ \wedge ~ \Xi (\hat{G}_{\rm{av}},e_{\rm{av}}) <0 \right\}\,, \label{eq:tk+1_event_av}
				\end{align}
				which is the average static event-trigger mechanism.
		\end{itemize}
		\item A feedback control action updated at the triggering instants:
		\begin{align}
			u_{k}=K\hat{G}_{\rm av}(\bar{t}_{k}) \,, \label{eq:U_MD4}
		\end{align}
		for all $\bar{t} \in \lbrack \bar{t}_{k}\,, \bar{t}_{k+1}\phantom{(}\!\!)$, $k\in \mathbb{N}$.
\end{enumerate}
\end{definition}

\begin{definition}[\small{Average Dynamic Triggering Condition}] \label{def:averageDynamicEvent}
Let $\Xi(\hat{G}_{\rm{av}},e_{\rm{av}})$ in (\ref{eq:Xi_event_2}) be the nonlinear mapping  and $K$  the control gain in (\ref{eq:U_event}), $\mu>0$  be a positive constant and $\upsilon_{\rm{av}}(\bar{t})$ be the solution of the dynamics
\begin{align}
\frac{d\upsilon_{\rm{av}}(\bar{t})}{d\bar{t}}&=-\frac{\mu}{\omega} \upsilon_{\rm{av}}(\bar{t}) +\frac{1}{\omega}\Xi (\hat{G}_{\rm{av}},e_{\rm{av}}) \,, ~ \mu>0\,, ~ \upsilon(0)\geq 0\,. \label{eq:dotUpsilon_ave}
\end{align}
Consequently, the event-triggered controller with dynamic-triggering condition consists of two components:
\begin{enumerate}
	\item A set of increasing sequence of time $I=\{\bar{t}_{0}\,, \bar{t}_{1}\,, \bar{t}_{2}\,,\ldots\}$ with $\bar{t}_{0}=0$ generated under the following rule:
		\begin{itemize}
			\item If $\left\{\bar{t} \in\mathbb{R}^{+}: \bar{t}\mathbb{>}\bar{t}_{k} \wedge \upsilon_{\rm{av}}(\bar{t})\mathbb{+}\gamma\Xi (\hat{G}_{\rm{av}},e_{\rm{av}}) \mathbb{<} 0\right\} \mathbb{=} \emptyset$, then the set of the times of the events is $I=\{\bar{t}_{0}\,, \bar{t}_{1}\,, \ldots, \bar{t}_{k}\}$.
			\item If $\left\{\bar{t} \in\mathbb{R}^{+}: \bar{t}\mathbb{>}\bar{t}_{k} \wedge \upsilon_{\rm{av}}(\bar{t})\mathbb{+}\gamma\Xi (\hat{G}_{\rm{av}},e_{\rm{av}}) \mathbb{<} 0\right\} \mathbb{\neq} \emptyset$, next event time is given by
				\begin{align}
					\!\!\!\!\!\!\!\! \bar{t}_{k+1}&\mathbb{=}\inf\left\{\bar{t} \in\mathbb{R}^{+}: \bar{t}\mathbb{>}\bar{t}_{k}  \wedge  \upsilon_{\rm{av}}(\bar{t})\mathbb{+}\gamma\Xi (\hat{G}_{\rm{av}},e_{\rm{av}})  \mathbb{<}  0 \right\}\!\!, \label{eq:tk+1_event_dynamic_av}
				\end{align}
				which is the average dynamic event-trigger mechanism.
		\end{itemize}
	\item A feedback control action updated at the triggering instants given by \eqref{eq:U_MD4}.
\end{enumerate}  
For all $\bar t\in(\bar t_{k},\bar t_{k+1})$, $\upsilon_{\rm{av}}(\bar t_{0})=\upsilon_{\rm{av}}(0)$ and $\upsilon_{\rm{av}}(\bar t_{k}^{-})=\upsilon_{\rm{av}}(\bar t_{k})=\upsilon_{\rm{av}}(\bar t_{k}^{+})$.
\end{definition}

We claim that the two event-triggering mechanisms discussed above  guarantee the asymptotic stabilization of $\hat{G}_{\text{av}}(\bar{t})$ and, consequently, that of $\tilde{\theta}_{\text{av}}(\bar{t})$. Since $H^{\ast}$ is invertible, both $\hat{G}_{\text{av}}(\bar{t})$ and $\tilde{\theta}_{\text{av}}(\bar{t})$ converge to the origin according to the averaging theory \cite{K:2002}.

\begin{remark}
From the time-scaling relation $\bar{t}=\omega t$, where $\omega$ is a constant, the static and dynamic triggering mechanisms of the averaging and the original system  are equivalent and only differ in the average sense.
\end{remark}

Next, the stability analysis will be carried out considering the static and dynamic event-triggering mechanisms. Note that, in both strategies it is considered the total lack of knowledge of the nonlinear map (\ref{eq:Q_1_event}), {\it i.e.}, the Hessian $H^{\ast}$, the optimizer $\theta^{\ast}$ and the extremum $Q^{\ast}$ are unknown parameters.

\textcolor{black}{In Sections~\ref{ETESNC_unknownH*} and \ref{DETESNC_unknownH*}, when analyzing the closed-loop average system, we employ and adapt the Lyapunov-based formalism originally presented
in the papers by Tabuada \cite{T:2007} and Girard \cite{G:2014}, for the static and dynamic triggering mechanisms, respectively. Indeed, the beauty of the proposed approach is that our average dynamics in both analyses gracefully match those of Tabuada and Girard.} 

\section{Static Event-Triggering in Extremum Seeking}\label{ETESNC_unknownH*}

\begin{figure}[H]
\includegraphics[width=8.5cm]{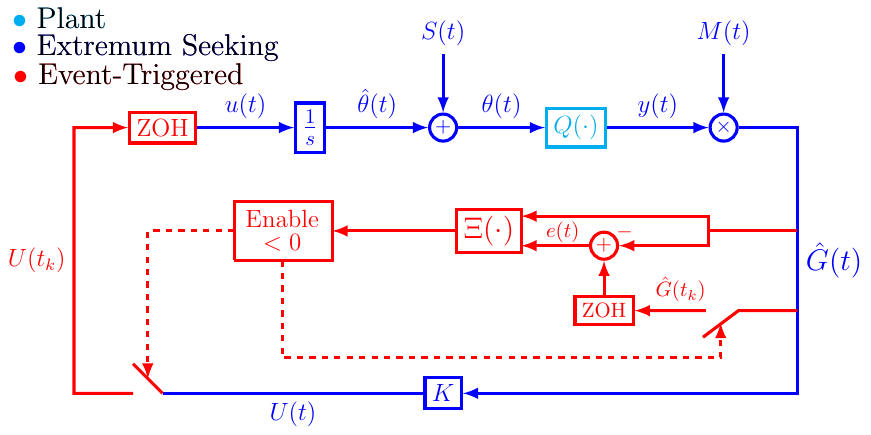}
\caption{Event-triggered based on extremum seeking scheme.}
\label{fig:blockDiagram_8}
\end{figure}

Theorem~\ref{thm:NETESC_2} states the local asymptotic stability of the extremum seeking based on dynamic event-triggered execution mechanism shown in Fig.~\ref{fig:blockDiagram_8} is ensured.
\textcolor{black}{
\begin{theorem} \label{thm:NETESC_2}
Consider the closed-loop average dynamics of the gradient estimate (\ref{eq:dotHatGav_event_1}), the average error vector \eqref{eq:Eav_event_1} and the average \textbf{static} event-triggered mechanism given by \textbf{Definition \ref{def:averageStaticEvent}}. Under  Assumptions \ref{assumption1}--\ref{assumption3} and considering the quadratic mapping $\Xi(\hat{G}_{\rm{av}},e_{\rm{av}})$ given by (\ref{eq:Xi_event_2}), for $\omega>0$, defined in (\ref{eq:omega_event_1}), sufficiently large, the equilibrium $\hat{G}_{\text{av}}(t)=0$ 
is locally exponentially stable and $\tilde{\theta}_{\text{av}}(t)$ converges exponentially to zero. In particular, there exist constants  $m\,,M_{\theta}\,,M_{y}>0$ such that
\begin{align}
\|\theta(t)-\theta^{\ast}\|&\leq M_{\theta}\exp(-mt)+\mathcal{O}\left(a+\frac{1}{\omega}\right)\,, \label{eq:normTheta_thm2} \\ 
|y(t)\mathbb{-}Q^{\ast}|&\leq M_{y}\exp(-mt)+\mathcal{O}\left(  a^{2} +\frac{1}{\omega^{2}}\right)\,, \label{eq:normY_thm2}
\end{align}
where $a=\sqrt{\sum_{i=1}^{n}a_{i}^{2}}$, with $a_i$ defined in \eqref{eq:S_event} and the constants $M_{\theta}$, and $M_{y}$ depending on the initial condition $\theta(0)$. In addition, there exists a lower bound  $\tau^{\ast}$ for the inter-execution interval $t_{k+1}-t_{k}$  for all $k \in \mathbb{N}$ precluding the Zeno behavior.
\end{theorem}
}
\textit{Proof:} The proof of the theorem is divided into two parts: stability analysis and avoidance of Zeno behavior.

\begin{flushleft}
\textcolor{black}{\underline{\it A. Stability Analysis}}
\end{flushleft}

Now, consider the following candidate Lyapunov function for the average system \eqref{eq:dotHatGav_event_1}
\begin{align}
V_{\text{av}}(\bar{t})=\hat{G}^{T}_{\text{av}}(\bar{t})P\hat{G}_{\text{av}}(\bar{t}) \,,\, P=P^T>0, \label{eq:Vav_event_pf2}
\end{align}  
with time-derivative
\begin{align}
\frac{dV_{\text{av}}(\bar{t})}{d\bar{t}}&=-\frac{1}{\omega}\hat{G}_{\text{av}}^{T}(\bar{t})Q\hat{G}_{\text{av}}(\bar{t})+\frac{1}{\omega}e_{\text{av}}^{T}(\bar{t})K^{T}H^{\ast T}P\hat{G}_{\text{av}}(\bar{t}) \nonumber \\
&\quad+\frac{1}{\omega}\hat{G}_{\text{av}}^{T}(\bar{t})PH^{\ast}Ke_{\text{av}}(\bar{t})\,, \label{eq:dotVav_event_1_pf2}
\end{align}
whose  upper bound satisfies
\begin{align}
\frac{dV_{\text{av}}(\bar{t})}{d\bar{t}}&\mathbb{\leq}\mathbb{-}\frac{\lambda_{\min}(Q)}{\omega}\|\hat{G}_{\text{av}}(\bar{t})\|^{2}\mathbb{+}
\frac{ 2\|PH^{\ast}K\|}{\omega}\|e_{\text{av}}(\bar{t})\| \|\hat{G}_{\text{av}}(\bar{t})\|. \label{eq:dotVav_event_2_pf2}
\end{align}
Using Assumptions \ref{assumption2} and \ref{assumption3}, we arrive at
\begin{align}
\frac{dV_{\text{av}}(\bar{t})}{d\bar{t}}&\leq-\frac{\alpha}{\omega}\|\hat{G}_{\text{av}}(\bar{t})\|^{2}+\frac{\beta}{\omega}\|e_{\text{av}}(\bar{t})\| \|\hat{G}_{\text{av}}(\bar{t})\| \,. \label{eq:dotVav_event_3_pf2}
\end{align}

In the proposed event-triggered mechanism, the update law is (\ref{eq:tk+1_event_av}) and $\Xi(\hat{G}_{\rm{av}},e_{\rm{av}})$ is given by (\ref{eq:Xi_event_2}). The signal $u_{\text{av}}(t)$ is held constant between two consecutive events, {\it i.e.}, while $\Xi(\hat{G}_{\rm{av}},e_{\rm{av}})\geq 0$, one has
\begin{align}
\Xi(\hat{G}_{\rm{av}},e_{\rm{av}}) &= \sigma\alpha\|\hat{G}_{\text{av}}(\bar{t})\|^{2}-\beta\|e_{\text{av}}(\bar{t})\| \|\hat{G}_{\text{av}}(\bar{t})\|\nonumber \\
&= \beta\|\hat{G}_{\text{av}}(\bar{t})\|\left(\frac{\sigma\alpha}{\beta}\|\hat{G}_{\text{av}}(\bar{t})\|-\|e_{\text{av}}(\bar{t})\| \right)\geq 0		\,.
\end{align}
Therefore, considering the event-triggered approach, the average measurement error $e_{\text{av}}(\bar{t})$ is upper bounded by
\begin{align}
\|e_{\text{av}}(\bar{t})\|&\leq \frac{\sigma \alpha}{\beta}\|\hat{G}_{\text{av}}(\bar{t})\|	\,. \label{eq:eAv_upperBound}
\end{align}
Now, plugging (\ref{eq:eAv_upperBound}) into (\ref{eq:dotVav_event_3_pf2}), 
\begin{align}
\frac{dV_{\text{av}}(\bar{t})}{d\bar{t}}&\leq-\frac{\alpha(1-\sigma)}{\omega}\|\hat{G}_{\text{av}}(\bar{t})\|^{2}\,. \label{eq:dotVav_event_4_pf2}
\end{align}

By using the Rayleigh-Ritz Inequality \cite{K:2002}, 
\begin{align}
\lambda_{\min}(P)\|\hat{G}_{\text{av}}(\bar{t})\|^{2}\leq V_{\text{av}}(\bar{t}) \leq \lambda_{\max}(P)\|\hat{G}_{\text{av}}(\bar{t})\|^{2}\,, \label{eq:Rayleigh-Ritz_pf2}
\end{align}
 and the following upper bound for (\ref{eq:dotVav_event_4_pf2})
\begin{align}
\frac{dV_{\text{av}}(\bar{t})}{d\bar{t}}&\leq -\frac{\alpha(1-\sigma)}{\omega}\|\hat{G}_{\text{av}}(\bar{t})\|^{2} \leq -\frac{\alpha(1-\sigma)}{\omega \lambda_{\max}(P)}V_{\text{av}}(\bar{t}) \,. \label{eq:dotVav_event_5_pf2}
\end{align}
Then, invoking the Comparison Principle \cite[Comparison Lemma]{K:2002} an upper bound $\bar{V}_{\text{av}}(\bar{t})$ for $V_{\text{av}}(\bar{t})$
\begin{align}
V_{\text{av}}(\bar{t})\leq \bar{V}_{\text{av}}(\bar{t}) \,, \quad \forall \bar{t}\in \lbrack \bar{t}_{k},\bar{t}_{k+1}\phantom{(}\!\!) \,, \label{eq:VavBarVav_1_pf2}
\end{align}
 is given by the solution of the equation
\begin{align}
\frac{d\bar{V}_{\text{av}}(\bar{t})}{d\bar{t}}=-\frac{\alpha(1-\sigma)}{\omega \lambda_{\max}(P)}\bar{V}_{\text{av}}(\bar{t})\,, \quad \bar{V}_{\text{av}}(\bar{t}_{k})=V_{\text{av}}(\bar{t}_{k})\,.
\end{align}
In other words, $ \forall \bar{t}\in \lbrack \bar{t}_{k},\bar{t}_{k+1}\phantom{(}\!\!)$,
\begin{align}
\bar{V}_{\text{av}}(\bar{t})=\exp\left(-\frac{\alpha(1-\sigma)}{\omega \lambda_{\max}(P)}\bar{t}\right)V_{\text{av}}(\bar{t}_{k})\,, \label{eq:_pf2}
\end{align}
and inequality (\ref{eq:VavBarVav_1_pf2}) is rewritten as
\begin{align}
V_{\text{av}}(\bar{t})\leq \exp\left(-\frac{\alpha(1-\sigma)}{\omega \lambda_{\max}(P)}\bar{t}\right)V_{\text{av}}(\bar{t}_{k}) \,, \quad \forall \bar{t}\in \lbrack \bar{t}_{k},\bar{t}_{k+1}\phantom{(}\!\!)
\,. \label{eq:VavBarVav_2_pf2}
\end{align}

By defining, $\bar{t}_{k}^{+}$ and $\bar{t}_{k}^{-}$ as the right and left limits of $\bar{t}=\bar{t}_{k}$, respectively, it easy to verify that $V_{\text{av}}(\bar{t}_{k+1}^{-})\leq \exp\left(-\frac{\alpha(1-\sigma)}{\omega \lambda_{\max}(P)}(\bar{t}_{k+1}^{-}-\bar{t}_{k}^{+})\right)V_{\text{av}}(\bar{t}_{k}^{+})$. Since $V_{\text{av}}(\bar{t})$ is continuous, $V_{\text{av}}(\bar{t}_{k+1}^{-})=V_{\text{av}}(\bar{t}_{k+1})$ and $V_{\text{av}}(\bar{t}_{k}^{+})=V_{\text{av}}(\bar{t}_{k})$, and therefore,
\begin{align}
    V_{\text{av}}(\bar{t}_{k+1})\leq \exp\left(-\frac{\alpha(1-\sigma)}{\omega \lambda_{\max}(P)}(\bar{t}_{k+1}-\bar{t}_{k})\right)V_{\text{av}}(\bar{t}_{k})\,. \label{METES_eq:mmd_1_s}
\end{align}
Hence, for any $\bar{t}\geq 0$ in $ \bar{t}\in \lbrack \bar{t}_{k},\bar{t}_{k+1}\phantom{(}\!\!)$, $k \in \mathbb{N}$, one has 
\begin{align}
    V_{\text{av}}(\bar{t})&\leq \exp\left(-\frac{\alpha(1-\sigma)}{\omega \lambda_{\max}(P)}(\bar{t}-\bar{t}_{k})\right) V_{\text{av}}(\bar{t}_{k}) \nonumber \\
    &\leq \exp\left(-\frac{\alpha(1-\sigma)}{\omega \lambda_{\max}(P)}(\bar{t}-\bar{t}_{k})\right) \nonumber \\
		&\quad \times \exp\left(-\frac{\alpha(1-\sigma)}{\omega \lambda_{\max}(P)}(\bar{t}_{k}-\bar{t}_{k-1})\right)V_{\text{av}}(\bar{t}_{k-1}) \nonumber \\
    &\leq \ldots \leq \nonumber \\
    &\leq \exp\left(-\frac{\alpha(1-\sigma)}{\omega \lambda_{\max}(P)}(\bar{t}\mathbb{-}\bar{t}_{k})\right) \nonumber \\
		&\quad \times \prod_{i=1}^{i=k}\exp\left(-\frac{\alpha(1-\sigma)}{\omega \lambda_{\max}(P)}(\bar{t}_{i}-\bar{t}_{i-1})\right)V_{\text{av}}(\bar{t}_{i-1}) \nonumber \\
    &=\exp\left(-\frac{\alpha(1-\sigma)}{\omega \lambda_{\max}(P)}\bar{t}\right) V_{\text{av}}(0)\,, \quad \forall \bar{t}\geq 0\,. \label{METES_eq:VavBarVav_2_pf2}
\end{align}

Now, lower bounding the left-hand side and upper bounding the right-hand size of (\ref{METES_eq:VavBarVav_2_pf2}) with the corresponding sides of (\ref{eq:Rayleigh-Ritz_pf2}), one gets
\begin{align}
\lambda_{\min}(P)\|\hat{G}_{\text{av}}(\bar{t})\|^{2}&\leq \exp\left(-\frac{\alpha(1-\sigma)}{\omega \lambda_{\max}(P)}\bar{t}\right) \nonumber \\
&\quad\times\lambda_{\max}(P)\|\hat{G}_{\text{av}}(0)\|^{2} \,. \label{eq:VavBarVav_3_pf2}
\end{align}
Then,
\begin{align}
&\|\hat{G}_{\text{av}}(\bar{t})\|^{2}\leq \exp\left(-\frac{\alpha(1-\sigma)}{\omega \lambda_{\max}(P)}\bar{t}\right)\frac{\lambda_{\max}(P)}{\lambda_{\min}(P)}\|\hat{G}_{\text{av}}(0)\|^{2} \nonumber \\
&=\left[\exp\left(-\frac{\alpha(1-\sigma)}{2\omega \lambda_{\max}(P)}\bar{t}\right)\sqrt{\frac{\lambda_{\max}(P)}{\lambda_{\min}(P)}}\|\hat{G}_{\text{av}}(0)\|\right]^{2}\,, \label{eq:VavBarVav_4_pf2}
\end{align}
and
\begin{align}
\|\hat{G}_{\text{av}}(\bar{t})\|\leq\exp\left(-\frac{\alpha(1-\sigma)}{2\omega\lambda_{\max}(P)}\bar{t}\right)\sqrt{\frac{\lambda_{\max}(P)}{\lambda_{\min}(P)}}\|\hat{G}_{\text{av}}(0)\|\,. \label{eq:normHatGav_1_pf2}
\end{align}
Although the analysis has been focused on the convergence of $\hat{G}_{\text{av}}(\bar{t})$ and, consequently, of $\hat{G}(t)$, the obtained results through (\ref{eq:normHatGav_1_pf2}) can be easily extended to the variables $\tilde{\theta}_{\text{av}}(\bar{t})$ and $\tilde{\theta}(t)$. Using relation (\ref{eq:hatGav_event_1}),   (\ref{eq:Vav_event_pf2}) and   (\ref{eq:dotVav_event_2_pf2}) are rewritten as
\begin{align}
V_{\text{av}}(\bar{t})&=
\tilde{\theta}^{T}_{\text{av}}(\bar{t})H^{\ast T}PH^{\ast}\tilde{\theta}_{\text{av}}(\bar{t})\,, \label{eq:Vav_event_2_pf2} \\
\frac{dV_{\text{av}}(\bar{t})}{d\bar{t}}&\leq-\frac{\alpha(1-\sigma)}{\omega \lambda_{\max}(P)}V_{\text{av}}(\bar{t}) \,. \label{eq:dotVav_event_7_pf2}
\end{align}
From Assumption \ref{assumption1}, the quadratic matrix $H^{\ast}$ has linearly independent rows and columns. Furthermore, from Assumption \ref{assumption2}, $P$ is a symmetric and positive definite matrix. Thus, there exists a matrix $R$ with independent columns such that $P=R^{T}R$ and, consequently, $\bar{P}=H^{\ast T}PH^{\ast}$ is a symmetric and positive definite matrix \cite[Section~6.5]{S:2016} leading to 
\begin{align}
V_{\text{av}}(\bar{t})&=\tilde{\theta}^{T}_{\text{av}}(\bar{t})\bar{P}\tilde{\theta}_{\text{av}}(\bar{t})\,, \label{eq:Vav_event_3_pf2}\,
\end{align}
satisfying the Rayleigh-Ritz inequality
\begin{align}
\lambda_{\min}(\bar{P})\|\tilde{\theta}_{\text{av}}(\bar{t})\|^2 \leq V_{\rm{av}} \leq \lambda_{\max}(\bar{P})\|\tilde{\theta}_{\text{av}}(\bar{t})\|^2 \,. \label{eq:RRI_tildeTheta_pf2}
\end{align}
Then, by using (\ref{eq:VavBarVav_2_pf2}), (\ref{eq:Vav_event_3_pf2}) and (\ref{eq:RRI_tildeTheta_pf2}), one has
\begin{align}
\lambda_{\min}(\bar{P})\|\tilde{\theta}_{\text{av}}(\bar{t})\|^2&\leq\exp\left(-\frac{(1-\sigma)\lambda_{\min}(Q)}{\omega\lambda_{\max}(P)}\bar{t}\right) \nonumber \\
&\quad \times \lambda_{\max}(\bar{P})\|\tilde{\theta}_{\text{av}}(0)\|^2\,, \label{eq:normTildeThetaAv_0_pf2}
\end{align}
and
\begin{align}
\|\tilde{\theta}_{\text{av}}(\bar{t})\|\leq\exp\left(-\frac{\alpha(1-\sigma)}{2\omega\lambda_{\max}(P)}\bar{t}\right)\sqrt{\frac{\lambda_{\max}(\bar{P})}{\lambda_{\min}(\bar{P})}}\|\tilde{\theta}_{\text{av}}(0)\|\,. \label{eq:normTildeThetaAv_1_pf2}
\end{align}

Since (\ref{eq:dotTildeTheta_3_event}) has discontinuous right-hand side and the mapping $\tilde{\Theta}(\bar{t},\tilde{\theta}(\bar{t}),e(\bar{t}))$ in (\ref{eq:tildeMathcalTheta}) is $T$-periodic in $t$. From (\ref{eq:normTildeThetaAv_1_pf2}),  $\tilde{\theta}_{\text{av}}(\bar{t})$ is asymptotically stable, by invoking  \cite[Theorem~2]{P:1979}, such that      
\begin{align}
\|\tilde{\theta}(t)-\tilde{\theta}_{\text{av}}(t)\|\leq\mathcal{O}\left(\frac{1}{\omega}\right)\,.
\end{align}
By using the Triangle inequality \cite{A:1957}, one has
\begin{align}
&\|\tilde{\theta}(t)\|\leq\|\tilde{\theta}_{\text{av}}(t)\|+\mathcal{O}\left(\frac{1}{\omega}\right)\nonumber \\
&\leq \exp\left(-\frac{\alpha(1-\sigma)}{2\lambda_{\max}(P)}t\right)\sqrt{\frac{\lambda_{\max}(\bar{P})}{\lambda_{\min}(\bar{P})}}\|\tilde{\theta}_{\text{av}}(0)\|+\mathcal{O}\left(\frac{1}{\omega}\right)\!. \label{eq:nomrTildeTheta_pf2}
\end{align}
Furthermore, 
\begin{align}
\|\hat{G}(t)-\hat{G}_{\text{av}}(t)\|\leq\mathcal{O}\left(\frac{1}{\omega}\right)\,,
\end{align}
and by using again the Triangle inequality \cite{A:1957}, such that
\begin{align}
&\|\hat{G}(t)\|\leq\|\hat{G}_{\text{av}}(t)\|+\mathcal{O}\left(\frac{1}{\omega}\right)\nonumber \\
&\leq \exp\left(-\frac{\alpha(1-\sigma)}{2\lambda_{\max}(P)}t\right)\sqrt{\frac{\lambda_{\max}(P)}{\lambda_{\min}(P)}}\|\hat{G}_{\text{av}}(0)\|+\mathcal{O}\left(\frac{1}{\omega}\right)\!. \label{eq:nomrTildeTheta_pf2_}
\end{align}
Now, from (\ref{eq:theta_2_event}), we can write
\begin{align}
\theta(t)-\theta^{\ast}=\tilde{\theta}(t)+S(t)\,, \label{eq:theta_3_event_pf2}
\end{align}
whose norm satisfies
\begin{align}
&\|\theta(t)-\theta^{\ast}\|=\|\tilde{\theta}(t)+S(t)\| \leq \|\tilde{\theta}(t)\|+\|S(t)\| \nonumber \\
&\leq \exp\left(-\frac{\alpha(1-\sigma)}{2\lambda_{\max}(P)}t\right)\sqrt{\frac{\lambda_{\max}(\bar{P})}{\lambda_{\min}(\bar{P})}}\|\theta(0) - \theta^{\ast}\|\nonumber \\
&\quad+\mathcal{O}\left(a+\frac{1}{\omega}\right). \label{eq:theta_4_event_pf2}
\end{align}
Defining the error variable $\tilde{y}(t)$ as 
\begin{align}
\tilde{y}(t):=y(t)-Q^{\ast}\,,
\end{align}
using the fact that  $y(t)=Q(\theta(t))$ where $Q(\theta(t))$ is defined in \eqref{eq:Q_1_event} and the Cauchy-Schwarz inequality \cite{S:2008}, we get
\begin{align}
|\tilde{y}(t)|&=|y(t)-Q^{\ast}| 
=|(\theta(t)-\theta^{\ast})^{T}H^{\ast}(\theta(t)-\theta^{\ast})| \nonumber \\
&\leq \|H^{\ast}\|\|\theta(t)-\theta^{\ast}\|^{2}\,, \label{tildeY_event_1_pf2}
\end{align}
and substituting  (\ref{eq:theta_4_event_pf2}) in \eqref{tildeY_event_1_pf2}, the following holds:
\begin{align}
&|\tilde{y}(t)|\mathbb{\leq} \|H^{\ast}\|\! \left[\! \exp\left(\!\!-\frac{\alpha(1-\sigma)}{2\lambda_{\max}(P)}t\!\right)\sqrt{\frac{\lambda_{\max}(\bar{P})}{\lambda_{\min}(\bar{P})}}\|\theta(0) - \theta^{\ast}\|\right.\nonumber \\
&\quad\left.+\mathcal{O}\left(a+\frac{1}{\omega}\right)\right]^{2} \nonumber \\
&= \|H^{\ast}\| \left[ \exp\left(-\frac{\alpha(1-\sigma)}{\lambda_{\max}(P)}t\right)\frac{\lambda_{\max}(\bar{P})}{\lambda_{\min}(\bar{P})}\|\theta(0) - \theta^{\ast}\|^{2}\right.\nonumber \\
&\quad+2\exp\left(-\frac{\alpha(1-\sigma)}{2\lambda_{\max}(P)}t\right)\sqrt{\frac{\lambda_{\max}(\bar{P})}{\lambda_{\min}(\bar{P})}}\|\theta(0) - \theta^{\ast}\|  \nonumber \\
&\quad \left.\times\mathcal{O}\left(a+\frac{1}{\omega}\right)+\mathcal{O}\left(a+\frac{1}{\omega}\right)^{2}\right] \,. \label{tildeY_event_2_pf2}
\end{align}
Since $\exp\left(-\frac{\alpha(1-\sigma)}{\lambda_{\max}(P)}t\right)\leq \exp\left(-\frac{\alpha(1-\sigma)}{2\lambda_{\max}(P)}t\right)$ and, according to \cite[Definition~10.1]{K:2002}, $\|H^{\ast}\|\mathcal{O}\left(a+\frac{1}{\omega}\right)^{2}$ is of order of magnitude $\mathcal{O}\left(a+\frac{1}{\omega}\right)^{2}$,  (\ref{tildeY_event_1_pf2}) is upper bounded by
\begin{align}
&|y(t)-Q^{\ast}| \leq  \exp\left(-\frac{\alpha(1-\sigma)}{2\lambda_{\max}(P)}t\right)\|H^{\ast}\| \nonumber \\
&\times\left[\frac{\lambda_{\max}(\bar{P})}{\lambda_{\min}(\bar{P})}\|\theta(0) - \theta^{\ast}\|+2\sqrt{\frac{\lambda_{\max}(\bar{P})}{\lambda_{\min}(\bar{P})}}\mathcal{O}\left(a+\frac{1}{\omega}\right)\right] \nonumber \\
&\times\|\theta(0) - \theta^{\ast}\|+\mathcal{O}\left(a^{2}+2\frac{a}{\omega}+\frac{1}{\omega^{2}}\right) \,. \label{tildeY_event_3_pf2}
\end{align}
Finally, once $a\,,\omega>0$, by using the Young inequality \cite{K:2002}, $\frac{a}{\omega}\leq \frac{a^{2}}{2}+\frac{1}{2\omega^{2}}$, and
\begin{align}
&|y(t)-Q^{\ast}| \leq  \exp\left(-\frac{\alpha(1-\sigma)}{2\lambda_{\max}(P)}t\right)\|H^{\ast}\| \nonumber \\
&\times\left[\frac{\lambda_{\max}(\bar{P})}{\lambda_{\min}(\bar{P})}\|\theta(0) - \theta^{\ast}\|+2\sqrt{\frac{\lambda_{\max}(\bar{P})}{\lambda_{\min}(\bar{P})}}\mathcal{O}\left(a+\frac{1}{\omega}\right)\right] \nonumber \\
&\quad \times\|\theta(0) - \theta^{\ast}\|+\mathcal{O}\left(a^{2}+\frac{1}{\omega^{2}}\right) \,, \label{tildeY_event_4_pf2}
\end{align}
since $\mathcal{O}\left(2a^{2}+2/\omega\right)=2\mathcal{O}\left(a^{2}+1/\omega^{2}\right)$ has an order of magnitude of $\mathcal{O}\left(a^{2}+1/\omega^{2}\right)$ \cite[Definition~10.1]{K:2002}.

Therefore, by defining the positive constants
\begin{align} 
m&=\frac{\alpha(1-\sigma)}{2\lambda_{\max}(P)} \,,\label{eq:m_event_1_pf2} \\
M_{\theta}&=\sqrt{\frac{\lambda_{\max}(\bar{P})}{\lambda_{\min}(\bar{P})}}\|\theta(0) - \theta^{\ast}\|\,,\label{eq:M_theta_event_1_pf2} \\
M_{y}&=\|H^{\ast}\|\frac{\lambda_{\max}(\bar{P})}{\lambda_{\min}(\bar{P})}\|\theta(0) - \theta^{\ast}\|^{2}\nonumber \\
&\quad+2\|H^{\ast}\|\sqrt{\frac{\lambda_{\max}(\bar{P})}{\lambda_{\min}(\bar{P})}}\|\theta(0) - \theta^{\ast}\|\mathcal{O}\left(a+\frac{1}{\omega}\right)\,,\label{eq:M_y_event_1_pf2}
\end{align}
inequalities (\ref{eq:theta_4_event_pf2}) and (\ref{tildeY_event_4_pf2}) satisfy (\ref{eq:normTheta_thm2}) and (\ref{eq:normY_thm2}), respectively.  

\begin{flushleft}
\textcolor{black}{\underline{\it B. Avoidance of Zeno Behavior}}
\end{flushleft}

Since  the average closed-loop system consists of (\ref{eq:dotHatGav_event_1}), with the event-triggered mechanism  (\ref{eq:Xi_event_2}),  (\ref{eq:tk+1_event_av}), and a control signal's update law (\ref{eq:eAv_upperBound}) satisfying $\|e_{\text{av}}(\bar{t})\|>\frac{\sigma \alpha}{\beta}\|\hat{G}_{\text{av}}(\bar{t})\|$, we conclude that  $dV_{\text{av}}(\bar{t})/d\bar{t}<0$ from (\ref{eq:dotVav_event_3_pf2}) and for all $\bar{t} \in \lbrack t_{k}\,, t_{k+1}\lbrack$. Thus, one can state that
\begin{align}
\sigma\alpha \|\hat{G}_{\text{av}}(\bar{t})\|^{2}-\beta\|e_{\text{av}}(\bar{t})\|\|\hat{G}_{\text{av}}(\bar{t})\|\geq 0\,, \label{ineq:interEvents_1_static}
\end{align}
and using the Peter-Paul inequality \cite{W:1971}, $cd\leq \frac{c^2}{2\epsilon}+\frac{\epsilon d^2}{2}$ for all $c,d,\epsilon>0$, with $c=\|e_{\rm{av}}(\bar{t})\|$, $d=\|\hat{G}_{\rm{av}}(\bar{t})\|$ and $\epsilon=\frac{\sigma\alpha}{\beta}$, the inequality (\ref{ineq:interEvents_1_static}) is lower bounded by
\begin{align}
&\sigma\alpha \|\hat{G}_{\text{av}}(\bar{t})\|^{2}-\beta\|e_{\text{av}}(\bar{t})\|\|\hat{G}_{\text{av}}(\bar{t})\|\geq \nonumber \\
&\sigma\alpha \|\hat{G}_{\text{av}}(\bar{t})\|^{2}-\beta\left(\frac{\sigma\alpha}{2\beta}\|\hat{G}_{\rm{av}}(\bar{t})\|^2+\frac{\beta}{2\sigma\alpha}\|e_{\rm{av}}(\bar{t})\|^2\right)\nonumber \\
&= q\|\hat{G}_{\rm{av}}(\bar{t})\|^{2}-p\|e_{\rm{av}}(\bar{t})\|^2\,,\label{ineq:interEvents_2_static_pf2}
\end{align}
where 
\begin{align}
q&=\frac{\sigma\alpha}{2} \quad \mbox{and} \quad p=\frac{\beta^2}{2\sigma\alpha}\,.\label{ineq:interEvents_3_static_pf2}
\end{align} 
In \cite{G:2014}, it is shown that a lower bound for the inter-execution interval is given by the time duration it takes for the function
\begin{align}
\phi(\bar{t})=\sqrt{\frac{p}{q}}\frac{\|e_{\rm{av}}(\bar{t})\|}{\|\hat{G}_{\rm{av}}(\bar{t})\|} \label{eq:phi_1_static_pf2}
\end{align}
to go from 0 to 1. The time-derivative of (\ref{eq:phi_1_static_pf2}) is 
\begin{align}
\frac{d\phi(\bar{t})}{d\bar{t}}&=\sqrt{\frac{p}{q}}\frac{1}{\|e_{\rm{av}}(\bar{t})\|\|\hat{G}_{\rm{av}}(\bar{t})\|}\left[e_{\rm{av}}^{T}(\bar{t})\frac{de_{\rm{av}}(\bar{t})}{d\bar{t}}\right.\nonumber \\
&\quad \left.-\hat{G}_{\rm{av}}^{T}(\bar{t})\frac{d\hat{G}_{\rm{av}}(\bar{t})}{d\bar{t}}\left(\frac{\|e_{\rm{av}}(\bar{t})\|}{\|\hat{G}_{\rm{av}}(\bar{t})\|}\right)^2\right]\,. \label{eq:dotPhi_1_static_pf2}
\end{align}
Now, plugging equations (\ref{eq:dotHatGav_event_1}) and (\ref{eq:Eav_event_1}) into (\ref{eq:dotPhi_1_static_pf2}), one arrives to 
\begin{align}
&\frac{d\phi(\bar{t})}{d\bar{t}}=\frac{1}{\omega}\sqrt{\frac{p}{q}}\frac{1}{\|e_{\rm{av}}(\bar{t})\|\|\hat{G}_{\rm{av}}(\bar{t})\|}\left\{-e_{\rm{av}}^{T}(\bar{t})H^{\ast}Ke_{\rm{av}}(\bar{t})\right.\nonumber \\
&\quad -e_{\rm{av}}^{T}(\bar{t})H^{\ast}K\hat{G}_{\rm{av}}(\bar{t})-\left[\hat{G}_{\rm{av}}^{T}(\bar{t})H^{\ast}K\hat{G}_{\rm{av}}(\bar{t})\right. \nonumber \\
&\quad\left.\left. +\hat{G}_{\rm{av}}^{T}(\bar{t})H^{\ast}Ke_{\rm{av}}(\bar{t})\right]\left(\frac{\|e_{\rm{av}}(\bar{t})\|}{\|\hat{G}_{\rm{av}}(\bar{t})\|}\right)^2\right\}\,. \label{eq:dotPhi_1_1_static_pf2}
\end{align}
Then, the following estimate holds:
\begin{align}
&\frac{d\phi(\bar{t})}{d\bar{t}}\leq\frac{1}{\omega}\sqrt{\frac{p}{q}}\frac{\|H^{\ast}K\|}{\|e_{\rm{av}}(\bar{t})\|\|\hat{G}_{\rm{av}}(\bar{t})\|}\left\{\|e_{\rm{av}}(\bar{t})\|^2\right.\nonumber \\
&\quad +\|e_{\rm{av}}(\bar{t})\|\|\hat{G}_{\rm{av}}(\bar{t})\|+\left[\|\hat{G}_{\rm{av}}(\bar{t})\|^2\right. \nonumber \\
&\quad\left.\left. +\|\hat{G}_{\rm{av}}(\bar{t})\|\|e_{\rm{av}}(\bar{t})\|\right]\left(\frac{\|e_{\rm{av}}(\bar{t})\|}{\|\hat{G}_{\rm{av}}(\bar{t})\|}\right)^2\right\} \nonumber\\
&=\frac{\|H^{\ast}K\|}{\omega}\sqrt{\frac{p}{q}}\left\{1+2\frac{\|e_{\rm{av}}(\bar{t})\|}{\|\hat{G}_{\rm{av}}(\bar{t})\|}+\frac{\|e_{\rm{av}}(\bar{t})\|^2}{\|\hat{G}_{\rm{av}}(\bar{t})\|^2}\right\}\,. \label{eq:dotPhi_2_static_pf2}
\end{align}

Hence,  using (\ref{eq:phi_1_static_pf2}), inequality (\ref{eq:dotPhi_2_static_pf2}) is rewritten as
\begin{align}
\omega\frac{d\phi(\bar{t})}{d\bar{t}}&\leq\|H^{\ast}K\|\sqrt{\frac{p}{q}}+2\|H^{\ast}K\|\phi(\bar{t})+\|H^{\ast}K\|\sqrt{\frac{q}{p}}\phi^{2}(\bar{t})\,. \label{eq:dotPhi_3_static_pf2}
\end{align}
From  the time-scaling $t =\frac{\bar{t}}{\omega} $, inequality (\ref{eq:dotPhi_2_static_pf2}) and invoking the Comparison Lemma \cite{K:2002}, a lower bound for the inter-execution time is found as
\begin{align}
\tau^{\ast}=\int_{0}^{1}\dfrac{1}{b_{0}+b_{1}\xi+b_{2}\xi^2}{\it d}\xi\,,\label{eq:tauAst_static_pf2}
\end{align}
with $b_{0}=\dfrac{\beta \|H^{\ast}K\| }{\alpha\sigma}$, $b_{1}=2\|H^{\ast}K\|$ and $b_{2}=\dfrac{\alpha\|H^{\ast}K\|\sigma}{\beta}$. Therefore, the Zeno behavior  is avoided. 

\textcolor{black}{While we do establish a minimum switching time to rule out a
Zeno behavior, the oscillations of high-frequency nature due to the perturbation or dither signals $S(t)$ and $M(t)$, which serve as excitation probing for the adaptive algorithm, will still occur according to the definition in (\ref{eq:S_event}) and (\ref{eq:M_event}).} \hfill $\square$

\section{Dynamic Event-Triggering in Extremum Seeking} \label{DETESNC_unknownH*}

\begin{figure}[h!]
\includegraphics[width=8.5cm]{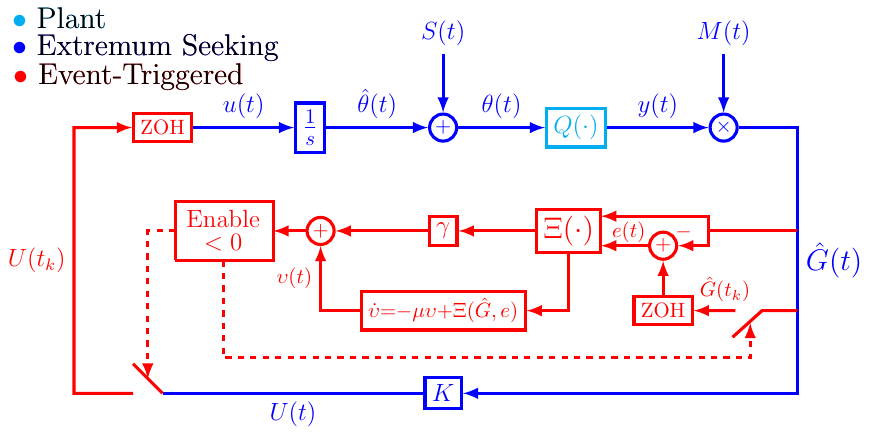}
\caption{Block diagram of the extremum seeking based on dynamic event-triggered mechanism.}
\label{fig:blockDiagram_10}
\end{figure}

\textcolor{black}{The motivation and need for the dynamic trigger is clear in the sense of reducing the number of switching and then not overloading the bandwidth of the channel to be used in the execution of the control loop. However, it is important to have in mind that this improvement may be followed by worse performances in terms of convergence rates and higher control efforts when compared to the previous the static-trigger approach.}

Theorem~\ref{thm:NETESC_4} demonstrates how  local asymptotic stability of the extremum seeking based on a dynamic event-triggered mechanism shown in Figure~\ref{fig:blockDiagram_10} is ensured.
\textcolor{black}{
\begin{theorem} \label{thm:NETESC_4}
Consider the closed-loop average dynamics of the gradient estimate (\ref{eq:dotHatGav_event_1}), the average error vector \eqref{eq:Eav_event_1} and the average \textbf{dynamic} event-triggered mechanism given by \textbf{Definition \ref{def:averageDynamicEvent}}. Under  Assumptions \ref{assumption1}--\ref{assumption3} and considering the quadratic mapping $\Xi(\hat{G}_{\rm{av}},e_{\rm{av}})$ given by (\ref{eq:Xi_event_2}), for  $\omega>0$,\linebreak defined in (\ref{eq:omega_event_1}), sufficiently large, 
the equilibrium $\hat{G}_{\text{av}}(t)=0$ is locally exponentially stable and  
$\tilde{\theta}_{\text{av}}(t)$ converges exponentially to zero. In particular, there exist constants  $m\,,M_{\theta}\,,M_{y}>0$ such that
\begin{align}
\|\theta(t)-\theta^{\ast}\|&\leq M_{\theta}\exp(-mt)+\mathcal{O}\left(a+\frac{1}{\omega}\right)\,, \label{eq:normTheta_thm4} \\ 
|y(t)\mathbb{-}Q^{\ast}|&\leq M_{y}\exp(-mt)+\mathcal{O}\left(  a^{2} +\frac{1}{\omega^{2}}\right)\,, \label{eq:normY_thm4}
\end{align}
where $a=\sqrt{\sum_{i=1}^{n}a_{i}^{2}}$, with $a_i$ defined in \eqref{eq:S_event} and the constants $M_{\theta}$ and $M_{y}$ depending on the initial condition $\theta(0)$. In addition, there exists a lower bound  $\tau^{\ast}$ for the inter-execution interval $t_{k+1}-t_{k}$  for all $k \in \mathbb{N}$ precluding the Zeno behavior.
\end{theorem}
}

\textit{Proof:} The proof of the theorem is again divided into two parts: stability analysis and avoidance of Zeno behavior.

\begin{flushleft}
\textcolor{black}{\underline{\it A. Stability Analysis}}
\end{flushleft}

First, notice that the dynamic triggering mechanism  (\ref{eq:tk+1_event_dynamic}) ensures, for all $t \in \lbrack t_{k}\,, t_{k+1}\lbrack$,
\begin{align}
\upsilon_{\rm{av}}(\bar{t})+\gamma \Xi(\hat{G}_{\rm{av}},e_{\rm{av}})\geq 0\,,
\end{align}
and,
\begin{align}
\Xi(\hat{G}_{\rm{av}},e_{\rm{av}})\geq -\frac{1}{\gamma}\upsilon_{\rm{av}}(\bar{t})\,. \label{eq:z_eta_1_pf4}
\end{align}
Now, with the help of  (\ref{eq:z_eta_1_pf4}), the following estimate of  (\ref{eq:dotUpsilon_ave}) holds
\begin{align}
\frac{d\upsilon_{\rm{av}}(\bar{t})}{d\bar{t}}&=-\frac{\mu}{\omega}\upsilon_{\rm{av}}(\bar{t})+\frac{1}{\omega}\Xi(\hat{G}_{\rm{av}},e_{\rm{av}}) \nonumber \\
&\geq-\frac{\mu}{\omega}\upsilon_{\rm{av}}(\bar{t})-\frac{1}{\omega\gamma}\upsilon_{\rm{av}}(\bar{t})=-\frac{1}{\omega}\left(\mu+\frac{1}{\gamma}\right)\upsilon_{\rm{av}}(\bar{t})\,. \label{ineq:dotEta_pf4}
\end{align}
Invoking \cite[Comparison Lemma, pp. 102]{K:2002}, the solution $\hat{\upsilon}_{\rm{av}}(\bar{t})$ of the following first-order dynamics
\begin{align}
\frac{d\hat{\upsilon}_{\rm{av}}(\bar{t})}{d\bar{t}}&=-\frac{1}{\omega}\left(\mu+\frac{1}{\gamma}\right)\hat{\upsilon}_{\rm{av}}(\bar{t})\,, \quad \hat{\upsilon}_{\rm{av}}(0)=\upsilon_{\rm{av}}(0)>0\,, \label{eq:dotHatEta_pf4}
\end{align}
precisely, 
\begin{align}
\hat{\upsilon}_{\rm{av}}(\bar{t})&=\exp\left(-\frac{1}{\omega}\left(\mu+\frac{1}{\gamma}\right)\bar{t}\right)\hat{\upsilon}_{\rm{av}}(0)>0\,, \quad \forall \bar{t}\geq 0\,, \label{eq:hatEta_pf4}
\end{align}
is a lower bound for $\upsilon_{\rm{av}}(\bar{t})$. To verify this fact, notice that, from (\ref{ineq:dotEta_pf4}) and (\ref{eq:dotHatEta_pf4}),
\begin{align}
\frac{d(\upsilon_{\rm{av}}(\bar{t})-\hat{\upsilon}_{\rm{av}}(\bar{t}))}{d\bar{t}}&\geq-\frac{1}{\omega}\left(\mu+\frac{1}{\gamma}\right)(\upsilon_{\rm{av}}(\bar{t})-\hat{\upsilon}_{\rm{av}}(\bar{t}))\,. \label{ineq:dotEta_dotHatEta_pf4}
\end{align}
Thus,
\begin{align}
\upsilon_{\rm{av}}(\bar{t})-\hat{\upsilon}_{\rm{av}}(\bar{t})&\geq \exp\left(-\frac{1}{\omega}\left(\mu+\frac{1}{\gamma}\right)\bar{t}\right)\underbrace{(\upsilon_{\rm{av}}(0)-\hat{\upsilon}_{\rm{av}}(0))}_{=0}
\end{align}
and 
\begin{align}
\upsilon_{\rm{av}}(\bar{t})&\geq \hat{\upsilon}_{\rm{av}}(\bar{t})>0\,, \quad \forall \bar{t}\geq 0\,. \label{eq:Eta_hatEta_pf4}
\end{align} 

Now, since $\upsilon_{\rm{av}}(\bar{t})>0$, for all $\upsilon_{\rm{av}}(\bar{t})\neq 0$, consider the following Lyapunov candidate for the   average system:
\begin{align}
V_{\rm{av}}(\bar{t})=\hat{G}^{T}_{\rm{av}}(\bar{t})P\hat{G}_{\rm{av}}(\bar{t})+\upsilon_{\rm{av}}(\bar{t})\,,\; P^{T}=P>0 \label{eq:lyapunov_dynamicETC_1_pf4}
\end{align}
The Rayleigh-Ritz inequality writes:
\begin{align}
\lambda_{\min}(P)\|\hat{G}_{\rm{av}}(\bar{t})\|^2 \leq \hat{G}^{T}_{\rm{av}}(\bar{t})P\hat{G}_{\rm{av}}(\bar{t}) \leq \lambda_{\max}(P)\|\hat{G}_{\rm{av}}(\bar{t})\|^2 \,. \label{ineq:RRI_dynamic_pf4}
\end{align}
The time-derivative of (\ref{eq:lyapunov_dynamicETC_1_pf4}) is given by
\begin{align}
&\frac{d{V}_{\text{av}}(\bar{t})}{d\bar{t}}=\frac{d\hat{G}_{\text{av}}^{T}(\bar{t})}{d\bar{t}}P\hat{G}_{\text{av}}(\bar{t})+\hat{G}^{T}_{\text{av}}(\bar{t})P\frac{d\hat{G}_{\text{av}}(\bar{t})}{d\bar{t}}+ \frac{d\upsilon_{\rm{av}}(\bar{t})}{d\bar{t}}\,,  \label{eq:dotLyapunov_dynamicETC_1_pf4}
\end{align} 
which, by using equations (\ref{eq:dotHatGav_event_1}) and (\ref{eq:dotUpsilon_ave}), can be rewritten as
\begin{align}
\frac{d{V}_{\text{av}}(\bar{t})}{d\bar{t}}&\mathbb{=}\mathbb{-}\frac{1}{\omega}\hat{G}_{\text{av}}^{T}(\bar{t})Q\hat{G}_{\text{av}}(\bar{t})\mathbb{+}\frac{1}{\omega}e_{\text{av}}^{T}(\bar{t})K^{T}H^{\ast T}P\hat{G}_{\text{av}}(\bar{t}) \nonumber \\
&\quad\mathbb{+}\frac{1}{\omega}\hat{G}_{\text{av}}^{T}(\bar{t})PH^{\ast}Ke_{\text{av}}(\bar{t}) \mathbb{-}\frac{\mu}{\omega}\upsilon_{\rm{av}}(\bar{t})\mathbb{+}\frac{1}{\omega}\Xi(\hat{G}_{\rm{av}},e_{\rm{av}})\,, \label{eq:dotLyapunov_dynamicETC_2_pf4}
\end{align}
Under Assumption~\ref{assumption2},  the following inequality is derived
\begin{align}
\frac{d{V}_{\text{av}}(\bar{t})}{d\bar{t}}
&\mathbb{\leq}\mathbb{-}\frac{\alpha}{\omega}\|\hat{G}_{\text{av}}(\bar{t})\|^2\mathbb{+}\frac{\beta}{\omega}\|e_{\text{av}}(\bar{t})\|\|\hat{G}_{\text{av}}(\bar{t})\|\mathbb{+}\frac{1}{\omega}\Xi (\hat{G}_{\rm{av}},e_{\rm{av}})\,. \label{eq:dotLyapunov_dynamicETC_2_2_pf4}
\end{align}
Plugging  (\ref{eq:Xi_event_2}) into (\ref{eq:dotLyapunov_dynamicETC_2_2_pf4}), one has
\begin{align}
\frac{d{V}_{\text{av}}(\bar{t})}{d\bar{t}}&\leq-\frac{(1-\sigma)\alpha}{\omega}\|\hat{G}_{\text{av}}(\bar{t})\|^2-\frac{\mu}{\omega}\upsilon_{\rm{av}}(\bar{t})\,. \label{eq:dotLyapunov_dynamicETC_3_pf4}
\end{align} 
Now,  using  (\ref{ineq:RRI_dynamic_pf4}) and  (\ref{eq:lyapunov_dynamicETC_1_pf4}), inequality (\ref{eq:dotLyapunov_dynamicETC_3_pf4}) can be upper bounded as follows
\begin{align}
&\frac{d{V}_{\text{av}}(\bar{t})}{d\bar{t}}\leq-\frac{(1-\sigma)\alpha}{\omega\lambda_{\max}(P)}\hat{G}_{\text{av}}^{T}(\bar{t})P\hat{G}_{\text{av}}(\bar{t})-\frac{\mu}{\omega}\upsilon_{\rm{av}}(\bar{t})\\
&\leq-\frac{1}{\omega}\min\left\{\frac{(1\!-\!\sigma)\alpha}{\lambda_{\max}(P)},\mu\right\}(\hat{G}_{\text{av}}^{T}(\bar{t})P\hat{G}_{\text{av}}(\bar{t})\!+\!\upsilon_{\rm{av}}(\bar{t})) \\
&\leq-\frac{1}{\omega}\min\left\{\frac{(1-\sigma)\alpha}{\lambda_{\max}(P)},\mu\right\}V_{\rm{av}}(\bar{t})\,, \quad \forall \bar{t}\in \lbrack \bar{t}_{k},\bar{t}_{k+1}\phantom{(}\!\!)\,. \label{eq:dotLyapunov_dynamicETC_3_pf4_}
\end{align}
Then, invoking the Comparison Principle \cite[Comparison Lemma]{K:2002}, an upper bound for $\bar{V}_{\rm{av}}(\bar{t})$ for $V_{\rm{av}}(\bar{t})$ 
\begin{align}
V_{\rm{av}}(\bar{t}) \leq \bar{V}_{\rm{av}}(\bar{t})\,, \quad \forall \bar{t}\in \lbrack \bar{t}_{k},\bar{t}_{k+1}\phantom{(}\!\!)\,, \label{ineq:VavBarVav_dynamic_pf4}
\end{align} 
is given by the solution of the dynamics
\begin{align}
\frac{d\bar{V}_{\rm{av}}(\bar{t})}{d\bar{t}}&\mathbb{=}-\frac{1}{\omega}\min\left\{\frac{(1-\sigma)\alpha}{\lambda_{\max}(P)},\mu\right\}\bar{V}_{\rm{av}}(\bar{t}),~\bar{V}_{\text{av}}(\bar{t}_{k})\mathbb{=}V_{\text{av}}(\bar{t}_{k}). \label{eq:barVav_0_dynamic_pf4}
\end{align}
That is, $ \forall \bar{t}\in \lbrack \bar{t}_{k},\bar{t}_{k+1}\phantom{(}\!\!)$, one has:
\begin{align}
\bar{V}_{\rm{av}}(\bar{t})=\exp\left(-\frac{1}{\omega}\min\left\{\frac{(1-\sigma)\alpha}{\lambda_{\max}(P)},\mu\right\}\bar{t}\right)\bar{V}_{\rm{av}}(\bar{t}_{k})\,. \label{eq:barVav_dynamic_pf4}
\end{align}
Using (\ref{eq:barVav_0_dynamic_pf4}) and (\ref{eq:barVav_dynamic_pf4}), the inequality (\ref{ineq:VavBarVav_dynamic_pf4}) is rewritten as 
\begin{align}
V_{\rm{av}}(\bar{t})\leq\exp\left(\!-\!\frac{1}{\omega}\min\left\{\frac{(1-\sigma)\alpha}{\lambda_{\max}(P)},\mu\right\}\bar{t}\right)V_{\rm{av}}(\bar{t}_{k})\,. \label{ineq:Vav_dynamic_pf4}
\end{align}
By defining, $\bar{t}_{k}^{+}$ and $\bar{t}_{k}^{-}$ as the right and left limits of $\bar{t}=\bar{t}_{k}$, respectively, it easy to verify that $V_{\text{av}}(\bar{t}_{k+1}^{-})\leq \exp\left(-\frac{1}{\omega}\min\left\{\frac{(1-\sigma)\alpha}{\lambda_{\max}(P)},\mu\right\}(\bar{t}_{k+1}^{-}-\bar{t}_{k}^{+})\right)V_{\text{av}}(\bar{t}_{k}^{+})$. Since $V_{\text{av}}(\bar{t})$ is continuous, $V_{\text{av}}(\bar{t}_{k+1}^{-})=V_{\text{av}}(\bar{t}_{k+1})$ and $V_{\text{av}}(\bar{t}_{k}^{+})=V_{\text{av}}(\bar{t}_{k})$, and therefore,
\begin{align}
    V_{\text{av}}(\bar{t}_{k+1})\mathbb{\leq} \exp\left(\!\!\mathbb{-}\frac{1}{\omega}\min\left\{\!\frac{(1\mathbb{-}\sigma)\alpha}{\lambda_{\max}(P)},\mu\!\right\}(\bar{t}_{k+1}\mathbb{-}\bar{t}_{k})\right)V_{\text{av}}(\bar{t}_{k}). \label{METES_eq:mmd_2_s}
\end{align}
Hence, for any $\bar{t}\geq 0$ in $ \bar{t}\in \lbrack \bar{t}_{k},\bar{t}_{k+1}\phantom{(}\!\!)$, $k \in \mathbb{N}$, one has 
\begin{align}
    &V_{\text{av}}(\bar{t})\leq \exp\left(-\frac{1}{\omega}\min\left\{\frac{(1-\sigma)\alpha}{\lambda_{\max}(P)},\mu\right\}(\bar{t}-\bar{t}_{k})\right) V_{\text{av}}(\bar{t}_{k}) \nonumber \\
    &\leq \exp\left(-\frac{1}{\omega}\min\left\{\frac{(1-\sigma)\alpha}{\lambda_{\max}(P)},\mu\right\}(\bar{t}-\bar{t}_{k})\right)  \nonumber \\
		& \times \exp\left(-\frac{1}{\omega}\min\left\{\frac{(1-\sigma)\alpha}{\lambda_{\max}(P)},\mu\right\}(\bar{t}_{k}-\bar{t}_{k-1})\right)V_{\text{av}}(\bar{t}_{k-1}) \nonumber \\
    &\leq \ldots \leq \nonumber \\
    &\leq \exp\left(-\frac{1}{\omega}\min\left\{\frac{(1-\sigma)\alpha}{\lambda_{\max}(P)},\mu\right\}(\bar{t}\mathbb{-}\bar{t}_{k})\right) \nonumber \\
		& \times \prod_{i=1}^{i=k}\exp\left(-\frac{1}{\omega}\min\left\{\frac{(1-\sigma)\alpha}{\lambda_{\max}(P)},\mu\right\}(\bar{t}_{i}-\bar{t}_{i-1})\right)V_{\text{av}}(\bar{t}_{i-1}) \nonumber \\
    &=\exp\left(-\frac{1}{\omega}\min\left\{\frac{(1-\sigma)\alpha}{\lambda_{\max}(P)},\mu\right\}\bar{t}\right) V_{\text{av}}(0)\,. \label{METES_eq:barVav_dynamic_pf4}
\end{align}

From (\ref{eq:lyapunov_dynamicETC_1_pf4}), it follows 
\begin{align}
\hat{G}^{T}_{\rm{av}}(\bar{t})P\hat{G}_{\rm{av}}(\bar{t})\leq V_{\rm{av}}(\bar{t})\,. \label{ineq:hatGav_Vav_dynamic_1_pf4}
\end{align}
Consequently, combining  (\ref{METES_eq:barVav_dynamic_pf4}) and (\ref{ineq:hatGav_Vav_dynamic_1_pf4}), one gets
\begin{align}
\hat{G}^{T}_{\rm{av}}(\bar{t})P\hat{G}_{\rm{av}}(\bar{t})&\leq \exp\left(\!-\!\frac{1}{\omega}\min\left\{\frac{(1-\sigma)\alpha}{\lambda_{\max}(P)},\mu\right\}\bar{t}\right)  \nonumber \\
&\quad \times V_{\rm{av}}(0)\,, \nonumber \\
&= \exp\left(\!-\!\frac{1}{\omega}\min\left\{\frac{(1-\sigma)\alpha}{\lambda_{\max}(P)},\mu\right\}\bar{t}\right)  \nonumber \\
& \times \left(\hat{G}^{T}_{\rm{av}}(0)P\hat{G}_{\rm{av}}(0)+\upsilon_{\rm{av}}(0)\right)\,. \label{ineq:hatGav_Vav_dynamic_2_pf4}
\end{align}
Since there exists a positive scalar $\kappa$ such that 
\begin{align}
\upsilon_{\rm{av}}(0)\leq \kappa \hat{G}^{T}_{\rm{av}}(0)P\hat{G}_{\rm{av}}(0)\,,
\end{align}
it is possible to write
\begin{align}
\hat{G}^{T}_{\rm{av}}(\bar{t})P\hat{G}_{\rm{av}}(\bar{t})&\leq \exp\left(\!-\!\frac{1}{\omega}\min\left\{\frac{(1-\sigma)\alpha}{\lambda_{\max}(P)},\mu\right\}\bar{t}\right)  \nonumber \\
& \times \left(1+\kappa\right)\hat{G}^{T}_{\rm{av}}(0)P\hat{G}_{\rm{av}}(0)\,, \label{ineq:hatGav_Vav_dynamic_2_pf4_}
\end{align} 
Therefore, from (\ref{ineq:RRI_dynamic_pf4}), one gets
\begin{align}
\lambda_{\min}(P)\|\hat{G}_{\rm{av}}(\bar{t})\|^2 & \!\leq \! \exp\left(\!\!-\!\frac{1}{\omega}\min\!\left\{\frac{(1\!-\!\sigma)\alpha}{\lambda_{\max}(P)},\mu\right\}\!\bar{t}\right)\!  \nonumber \\
&\times \left(1+\kappa\right)\lambda_{\max}(P)\|\hat{G}_{\rm{av}}(0)\|^2\,. \label{ineq:hatGav_Vav_dynamic_3_pf4}
\end{align}
Then,
\begin{align}
\|\hat{G}_{\rm{av}}(\bar{t})\|^2 &\leq \exp\left(\!-\!\frac{1}{\omega}\min\left\{\frac{(1-\sigma)\alpha}{\lambda_{\max}(P)},\mu\right\}\bar{t}\right)  \nonumber \\
&\quad \times \frac{\left(1+\kappa\right)\lambda_{\max}(P)}{\lambda_{\min}(P)}\|\hat{G}_{\rm{av}}(0)\|^2 \nonumber \\
&= \left[\exp\left(\!-\!\frac{1}{2\omega}\min\left\{\frac{(1-\sigma)\alpha}{\lambda_{\max}(P)},\mu\right\}\bar{t}\right)  \right. \nonumber \\
&\quad \left. \times \sqrt{\frac{\left(1+\kappa\right)\lambda_{\max}(P)}{\lambda_{\min}(P)}}\|\hat{G}_{\rm{av}}(0)\|\right]^2\,. \label{ineq:hatGav_Vav_dynamic_4_pf4}
\end{align}
equivalently,
\begin{align}
&\|\hat{G}_{\rm{av}}(\bar{t})\|^2-\left[\exp\left(\!-\!\frac{1}{2\omega}\min\left\{\frac{(1-\sigma)\alpha}{\lambda_{\max}(P)},\mu\right\}\bar{t}\right)  \right. \nonumber \\
&\qquad \left. \times \sqrt{\frac{\left(1+\kappa\right)\lambda_{\max}(P)}{\lambda_{\min}(P)}}\|\hat{G}_{\rm{av}}(0)\|\right]^2\leq0\,. \label{ineq:hatGav_Vav_dynamic_5_pf4}
\end{align}
Hence,
\begin{align}
&\quad\left[\|\hat{G}_{\rm{av}}(\bar{t})\|+\exp\left(\!-\!\frac{1}{2\omega}\min\left\{\frac{(1-\sigma)\alpha}{\lambda_{\max}(P)},\mu\right\}\bar{t}\right)  \right. \nonumber \\
&\qquad \left. \times \sqrt{\frac{\left(1+\kappa\right)\lambda_{\max}(P)}{\lambda_{\min}(P)}}\|\hat{G}_{\rm{av}}(0)\|\right] \nonumber \\
&\qquad \times\left[\|\hat{G}_{\rm{av}}(\bar{t})\|-\exp\left(\!-\!\frac{1}{2\omega}\min\left\{\frac{(1-\sigma)\alpha}{\lambda_{\max}(P)},\mu\right\}\bar{t}\right) \right. \nonumber \\
&\qquad\left. \times \sqrt{\frac{\left(1+\kappa\right)\lambda_{\max}(P)}{\lambda_{\min}(P)}}\|\hat{G}_{\rm{av}}(0)\|\right] \leq0 \label{ineq:hatGav_Vav_dynamic_6_pf4}
\end{align}
and
\begin{align}
\|\hat{G}_{\rm{av}}(\bar{t})\|&\leq\exp\left(\!-\!\frac{1}{2\omega}\min\left\{\frac{(1-\sigma)\alpha}{\lambda_{\max}(P)},\mu\right\}\bar{t}\right)\  \nonumber \\
&\quad \times \sqrt{\frac{\left(1+\kappa\right)\lambda_{\max}(P)}{\lambda_{\min}(P)}}\|\hat{G}_{\rm{av}}(0)\| \,. \label{ineq:hatGav_Vav_dynamic_7_pf4}
\end{align}
Although the analysis has been focused on the convergence of $\hat{G}_{\text{av}}(\bar{t})$ and, consequently, $\hat{G}(t)$, the obtained results through (\ref{ineq:hatGav_Vav_dynamic_7_pf4}) can be easily extended to the variables $\tilde{\theta}_{\text{av}}(\bar{t})$ and $\tilde{\theta}(t)$. From Assumption~\ref{assumption1}, the quadratic matrix $H^{\ast}$ has linearly independent rows and columns. Furthermore, from Assumption~\ref{assumption2}, $P$ is a symmetric and positive definite matrix. Thus, there exist a matrix $R$ with independent columns such that $P=R^{T}R$ and, consequently, $\bar{P}=H^{\ast T}PH^{\ast}$ is a symmetric and positive definite matrix \cite[Section~6.5]{S:2016}. Thus, by using (\ref{eq:hatGav_event_1}), the quadratic term $\hat{G}_{\text{av}}^{T}(\bar{t})P\hat{G}_{\text{av}}(\bar{t})$ in (\ref{eq:lyapunov_dynamicETC_1_pf4}) is written as
\begin{align} 
\hat{G}_{\text{av}}^{T}(\bar{t})P\hat{G}_{\text{av}}(\bar{t})
&=\tilde{\theta}_{\text{av}}^{T}(\bar{t})\bar{P}\tilde{\theta}_{\text{av}}(\bar{t})\,,
\end{align}
with the Rayleigh-Ritz inequality
\begin{align}
\lambda_{\min}(\bar{P})\|\tilde{\theta}_{\rm{av}}(\bar{t})\|^2 \leq \tilde{\theta}^{T}_{\rm{av}}(\bar{t})\bar{P}\tilde{\theta}_{\rm{av}}(\bar{t}) \leq \lambda_{\max}(\bar{P})\|\tilde{\theta}_{\rm{av}}(\bar{t})\|^2 \,. \label{ineq:RRI_tildeTheta_dynamic_pf4}
\end{align}

Therefore, inequality (\ref{ineq:hatGav_Vav_dynamic_2_pf4_}) can be rewritten as
\begin{align}
\tilde{\theta}^{T}_{\rm{av}}(\bar{t})\bar{P}\tilde{\theta}_{\rm{av}}(\bar{t})&\leq \exp\left(\!-\!\frac{1}{\omega}\min\left\{\frac{(1-\sigma)\alpha}{\lambda_{\max}(P)},\mu\right\}\bar{t}\right) \nonumber \\
& \times \left(1+\kappa\right)\tilde{\theta}^{T}_{\rm{av}}(0)\bar{P}\tilde{\theta}_{\rm{av}}(0)\,. \label{ineq:tildeTheta_dynamic_1_pf4}
\end{align}

Now, by considering inequalities (\ref{ineq:tildeTheta_dynamic_1_pf4}) and (\ref{ineq:RRI_tildeTheta_dynamic_pf4}), and following the steps between (\ref{ineq:hatGav_Vav_dynamic_2_pf4_}) and (\ref{ineq:hatGav_Vav_dynamic_7_pf4}), one obtains
\begin{align}
\|\tilde{\theta}_{\rm{av}}(\bar{t})\|&\leq\exp\left(\!-\!\frac{1}{2\omega}\min\left\{\frac{(1-\sigma)\alpha}{\lambda_{\max}(P)},\mu\right\}\bar{t}\right)  \nonumber \\
&\quad \times \sqrt{\frac{\left(1+\kappa\right)\lambda_{\max}(\bar{P})}{\lambda_{\min}(\bar{P})}}\|\tilde{\theta}_{\rm{av}}(0)\| \,. \label{ineq:tildeTheta_dynamic_2_pf4}
\end{align}

Since the differential equation (\ref{eq:dotTildeTheta_3_event}) has discontinuous right-hand size, $\tilde{\Theta}(\bar{t},\tilde{\theta}(\bar{t}),e(\bar{t}))$ in (\ref{eq:tildeMathcalTheta}) is $T$-periodic in $t$ and satisfy the Lipschitz condition. From (\ref{ineq:tildeTheta_dynamic_2_pf4}), by invoking  \cite[Theorem~2]{P:1979}, the $\tilde{\theta}_{\text{av}}(\bar{t})$ is asymptotically stable. That is       
\begin{align}
\|\tilde{\theta}(t)-\tilde{\theta}_{\text{av}}(t)\|\leq\mathcal{O}\left(\frac{1}{\omega}\right)\,.
\end{align}
Using the triangle inequality \cite{A:1957}, one has
\begin{align}
\|\tilde{\theta}(t)\|&\leq\|\tilde{\theta}_{\text{av}}(t)\|+\mathcal{O}\left(\frac{1}{\omega}\right)\nonumber \\
&\leq \exp\left(\!-\!\frac{1}{2\omega}\min\left\{\frac{(1-\sigma)\alpha}{\lambda_{\max}(P)},\mu\right\}\bar{t}\right) \nonumber \\
&\quad \times\sqrt{\frac{\left(1+\kappa\right)\lambda_{\max}(\bar{P})}{\lambda_{\min}(\bar{P})}}\|\tilde{\theta}_{\text{av}}(0)\|+\mathcal{O}\left(\frac{1}{\omega}\right)\!. \label{eq:nomrTildeTheta_pf4}
\end{align}
Furthermore, 
\begin{align}
\|\hat{G}(t)-\hat{G}_{\text{av}}(t)\|\leq\mathcal{O}\left(\frac{1}{\omega}\right)\,,
\end{align}
and again using the triangle inequality \cite{A:1957}, one obtains
\begin{align}
\|\hat{G}(t)\|&\leq\|\hat{G}_{\text{av}}(t)\|+\mathcal{O}\left(\frac{1}{\omega}\right)\nonumber \\
&\leq \exp\left(\!-\!\frac{1}{2\omega}\min\left\{\frac{(1-\sigma)\alpha}{\lambda_{\max}(P)},\mu\right\}\bar{t}\right) \nonumber \\
&\quad \times\sqrt{\frac{\left(1+\kappa\right)\lambda_{\max}(\bar{P})}{\lambda_{\min}(\bar{P})}}\|\tilde{\theta}_{\text{av}}(0)\|+\mathcal{O}\left(\frac{1}{\omega}\right)\!. \label{eq:nomrTildeTheta_pf4_}
\end{align}
Now, from (\ref{eq:theta_2_event}), we have
\begin{align}
\theta(t)-\theta^{\ast}=\tilde{\theta}(t)+S(t)\,, \label{eq:theta_3_event_pf4}
\end{align}
whose  norm satisfies 
\begin{align}
\|\theta(t)-\theta^{\ast}\|&=\|\tilde{\theta}(t)+S(t)\| \leq \|\tilde{\theta}(t)\|+\|S(t)\| \nonumber \\
&\leq \exp\left(\!-\!\frac{1}{2\omega}\min\left\{\frac{(1-\sigma)\alpha}{\lambda_{\max}(P)},\mu\right\}\bar{t}\right) \nonumber \\
&\quad\times\sqrt{\frac{\left(1+\kappa\right)\lambda_{\max}(\bar{P})}{\lambda_{\min}(\bar{P})}}\|\theta(0) - \theta^{\ast}\|\nonumber \\
&\quad+\mathcal{O}\left(a+\frac{1}{\omega}\right). \label{eq:theta_4_event_pf4}
\end{align}
Defining the error variable
\begin{align}
\tilde{y}(t):=y(t)-Q^{\ast}\,,\quad  y(t)=Q(\theta(t)),
\end{align}
and using  Cauchy-Schwartz inequality \cite{S:2008}, we get
\begin{align}
|\tilde{y}(t)|&=|y(t)-Q^{\ast}| =|(\theta(t)-\theta^{\ast})^{T}H^{\ast}(\theta(t)-\theta^{\ast})| \nonumber \\
&\leq \|H^{\ast}\|\|\theta(t)-\theta^{\ast}\|^{2}\,, \label{tildeY_event_1_pf4}
\end{align}
and  with the help of (\ref{eq:theta_4_event_pf4}) 
\begin{align}
|\tilde{y}(t)|&\leq \|H^{\ast}\| \left[ \exp\left(\!-\!\frac{1}{\omega}\min\left\{\frac{(1-\sigma)\alpha}{\lambda_{\max}(P)},\mu\right\}\bar{t}\right) \right.\nonumber \\
&\quad \times \frac{\left(1+\kappa\right)\lambda_{\max}(\bar{P})}{\lambda_{\min}(\bar{P})}\|\theta(0) - \theta^{\ast}\|^2 \nonumber \\
&\quad+2\exp\left(\!-\!\frac{1}{2\omega}\min\left\{\frac{(1-\sigma)\alpha}{\lambda_{\max}(P)},\mu\right\}\bar{t}\right) \nonumber \\
&\quad \times\sqrt{\frac{\left(1+\kappa\right)\lambda_{\max}(\bar{P})}{\lambda_{\min}(\bar{P})}}\|\theta(0) - \theta^{\ast}\|\mathcal{O}\left(a+\frac{1}{\omega}\right)\nonumber \\
&\quad\left.+\mathcal{O}\left(a+\frac{1}{\omega}\right)^{2}\right] \,. \label{tildeY_event_2_pf4}
\end{align}
Since $\exp\left(\!\!\mathbb{-}\!\!\min\!\left\{\!\!\frac{(1-\sigma)\alpha}{\lambda_{\max}(P)},\mu\!\right\}\!\!\frac{\bar{t}}{\omega}\right)\!\mathbb{\leq}\! \exp\!\left(\!\!\mathbb{-}\!\min\!\left\{\!\!\frac{(1-\sigma)\alpha}{\lambda_{\max}(P)},\mu\!\right\}\frac{\bar{t}}{2\omega}\right)$ and from \cite[Definition~10.1]{K:2002}, $\|H^{\ast}\|\mathcal{O}\left(a+\frac{1}{\omega}\right)^{2}$ is of order $\mathcal{O}\left(a+\frac{1}{\omega}\right)^{2}$. Hence,  (\ref{tildeY_event_1_pf4}) leads to
\begin{align}
|y(t)-Q^{\ast}| &\leq  \exp\left(\!\mathbb{-}\!\frac{1}{2\omega}\min\left\{\frac{(1-\sigma)\alpha}{\lambda_{\max}(P)},\mu\right\}\bar{t}\right) \nonumber \\
&\times\|H^{\ast}\|\left[\frac{\left(1+\kappa\right)\lambda_{\max}(\bar{P})}{\lambda_{\min}(\bar{P})}\|\theta(0) - \theta^{\ast}\|\right.\nonumber \\
& \left.+2\sqrt{\frac{\left(1+\kappa\right)\lambda_{\max}(\bar{P})}{\lambda_{\min}(\bar{P})}}\left(a+\frac{1}{\omega}\right)\right]\|\theta(0) - \theta^{\ast}\|\nonumber \\
&+\mathcal{O}\left(a^{2}+2\frac{a}{\omega}+\frac{1}{\omega^{2}}\right)\,. \label{tildeY_event_3_pf4}
\end{align}
Knowing that $a\,,\omega>0$ and using  Young's inequality $\frac{a}{\omega}\leq \frac{a^{2}}{2}+\frac{1}{2\omega^{2}}$, one gets
\begin{align}
|y(t)-Q^{\ast}| &\leq  \exp\left(\!-\!\frac{1}{2\omega}\min\left\{\frac{(1-\sigma)\alpha}{\lambda_{\max}(P)},\mu\right\}\bar{t}\right)\|H^{\ast}\| \nonumber \\
&\times\left[\frac{\left(1+\kappa\right)\lambda_{\max}(\bar{P})}{\lambda_{\min}(\bar{P})}\|\theta(0) - \theta^{\ast}\|\right.\nonumber \\
& \left.+2\sqrt{\frac{\left(1+\kappa\right)\lambda_{\max}(\bar{P})}{\lambda_{\min}(\bar{P})}}\left(a+\frac{1}{\omega}\right)\right]\|\theta(0) - \theta^{\ast}\|\nonumber \\
&+\mathcal{O}\left(a^{2}+\frac{1}{\omega^{2}}\right) \, \label{tildeY_event_4_pf4}
\end{align}
since $\mathcal{O}\left(2a^{2}+2/\omega\right)=2\mathcal{O}\left(a^{2}+1/\omega^{2}\right)$ is of  order $\mathcal{O}\left(a^{2}+1/\omega^{2}\right)$ \cite[Definition~10.1]{K:2002}.

Now, defining the positive constants:
\begin{align} 
m&=\frac{1}{2}\min\left\{\frac{(1-\sigma)\alpha}{\lambda_{\max}(P)},\mu\right\} \,,\label{eq:m_event_1_pf4} \\
M_{\theta}&=\sqrt{\frac{(1+\kappa)\lambda_{\max}(\bar{P})}{\lambda_{\min}(\bar{P})}}\|\theta(0) - \theta^{\ast}\|\,,\label{eq:M_theta_event_1_pf4} \\
M_{y}&=\|H^{\ast}\|\frac{\left(1+\kappa\right)\lambda_{\max}(\bar{P})}{\lambda_{\min}(\bar{P})}\|\theta(0) - \theta^{\ast}\|^2\nonumber \\
& \quad+2\|H^{\ast}\|\sqrt{\frac{\left(1+\kappa\right)\lambda_{\max}(\bar{P})}{\lambda_{\min}(\bar{P})}}\left(a+\frac{1}{\omega}\right)\|\theta(0) - \theta^{\ast}\|\,,\label{eq:M_y_event_1_pf4}
\end{align}
inequalities (\ref{eq:theta_4_event_pf4}) and (\ref{tildeY_event_4_pf4}) satisfy (\ref{eq:normTheta_thm4}) and (\ref{eq:normY_thm4}), respectively.

\begin{flushleft}
\textcolor{black}{\underline{\it B. Avoidance of Zeno Behavior}}
\end{flushleft}

Notice that, from (\ref{eq:Xi_event_2}) and (\ref{eq:tk+1_event_dynamic_av}), and using the Peter-Paul inequality \cite{W:1971}, we can write $cd\leq \frac{c^2}{2\epsilon}+\frac{\epsilon d^2}{2}$, for all $c,d,\epsilon>0$, with $c=\|e_{\rm{av}}(\bar{t})\|$, $d=\|\hat{G}_{\rm{av}}(\bar{t})\|$, $\epsilon=\frac{\sigma\lambda_{\min}(Q)}{2\|PH^{\ast}K\|}$ and $\bar{t}\in \lbrack t_{k}\,,t_{k+1}\phantom{(}\!\!\!)$. The following holds
\begin{align}
&\upsilon_{\rm{av}}(\bar{t})\mathbb{+}\gamma\left[\sigma \alpha \|\hat{G}_{\text{av}}(\bar{t})\|^2\mathbb{-}\beta \|e_{\text{av}}(\bar{t})\|\|\hat{G}_{\text{av}}(\bar{t})\|\right] \geq \nonumber \\
& \upsilon_{\rm{av}}(\bar{t})\mathbb{+}\gamma\left[\sigma\alpha\|\hat{G}_{\rm{av}}(\bar{t})\|^{2}\mathbb{-}\frac{\beta}{2}\!\left(\!\frac{\sigma\alpha}{\beta}\|\hat{G}_{\rm{av}}(\bar{t})\|^2\mathbb{+}\frac{\beta}{\sigma\alpha}\|e_{\rm{av}}(\bar{t})\|^2\!\right)\right] \nonumber \\
&= \upsilon_{\rm{av}}(\bar{t})+\gamma\left(q\|\hat{G}_{\rm{av}}(\bar{t})\|^{2}-p\|e_{\rm{av}}(\bar{t})\|^2\right)\,,\label{ineq:interEvents_2_dynamic_pf4}
\end{align}
where 
\begin{align}
q&=\frac{\sigma\alpha}{2} \quad \mbox{and} \quad p=\frac{\beta^2}{2\sigma\alpha}\,.\label{ineq:interEvents_3_dynamic_pf4}
\end{align} 

The minimum dwell-time of the event-triggered framework is given by the time it takes for the function
\begin{align}
\phi(\bar{t})=\frac{\sqrt{\gamma p}\|e_{\rm{av}}(\bar{t})\|}{\sqrt{\upsilon_{\rm{av}}(\bar{t})+\gamma q\|\hat{G}_{\rm{av}}(\bar{t})\|^{2}}}\,, \label{eq:phi_1_dynamic_pf4}
\end{align}
to go from 0 to 1. The derivative of $\phi(\bar{t})$ in (\ref{eq:phi_1_dynamic_pf4}) is given by
\begin{align}
\frac{d\phi(\bar{t})}{d\bar{t}}&=\frac{\sqrt{\gamma p}e_{\rm{av}}^{T}(\bar{t})\dfrac{d e_{\rm{av}}(\bar{t})}{d\bar{t}}}{\|e_{\rm{av}}(\bar{t})\|\sqrt{\upsilon_{\rm{av}}(\bar{t})+\gamma q\|\hat{G}_{\rm{av}}(\bar{t})\|^{2}}}\nonumber\\
&\quad-\frac{\sqrt{\gamma p}\|e_{\rm{av}}(\bar{t})\|}{2(\upsilon_{\rm{av}}(\bar{t})+\gamma q\|\hat{G}_{\rm{av}}(\bar{t})\|^{2})^{3/2}}\nonumber \\
&\quad\times\left(\frac{d\upsilon_{\rm{av}}(\bar{t})}{d\bar{t}}+\gamma q\hat{G}_{\rm{av}}^{T}(\bar{t})\frac{d\hat{G}_{\rm{av}}(\bar{t})}{d\bar{t}}\right)\,. \label{eq:dotPhi_1_dynamic_pf4}
\end{align}
Now, from  (\ref{eq:eta}), (\ref{eq:dotHatGav_event_1}), (\ref{eq:Eav_event_1}) and (\ref{ineq:interEvents_2_dynamic_pf4}), one  arrives at \begin{align}
\dfrac{d e_{\rm{av}}(\bar{t})}{d\bar{t}}&=-\dfrac{d \hat{G}_{\rm{av}}(\bar{t})}{d\bar{t}}\,,\\ 
\left\|\dfrac{d \hat{G}_{\rm{av}}(\bar{t})}{d\bar{t}}\right\|&\leq \dfrac{1}{\omega}\|H^{\ast}K\|\|\hat{G}_{\rm{av}}(\bar{t})\|+\dfrac{1}{\omega}\|H^{\ast}K\|\|e_{\rm{av}}(\bar{t})\|\,, \\
\dfrac{d\upsilon_{\rm{av}}(\bar{t})}{d\bar{t}}&\geq-\dfrac{\mu}{\omega}\upsilon_{\rm{av}}(\bar{t})+\dfrac{q}{\omega}\|\hat{G}_{\rm{av}}(\bar{t})\|^{2}-\dfrac{p}{\omega}\|e_{\rm{av}}(\bar{t})\|^2\, 
\end{align}
and from (\ref{eq:dotPhi_1_dynamic_pf4}) the following inequality holds
\begin{align}
&\frac{d\phi(\bar{t})}{d\bar{t}}\leq\frac{\|H^{\ast}K\|}{\omega}\sqrt{\frac{p}{q}}+\frac{\|H^{\ast}K\|}{\omega}\phi(\bar{t})+\frac{1}{2\omega\gamma}\phi^{3}(\bar{t})\nonumber\\
&+\frac{\|H^{\ast}K\|}{\omega}\sqrt{\frac{q}{p}}\phi^{2}(\bar{t})+\frac{\mu}{2\omega}\phi(\bar{t})\nonumber \\
&+\frac{\gamma q\|\hat{G}_{\rm{av}}(\bar{t})\|^{2}}{2\omega(\upsilon_{\rm{av}}(\bar{t})+\gamma q\|\hat{G}_{\rm{av}}(\bar{t})\|^{2})}\left(-\mu-\frac{1}{\gamma}+2\|H^{\ast}K\|\right)\phi(\bar{t})\,. \label{eq:dotPhi_2_dynamic_pf4}
\end{align}
Hence, if $\|H^{\ast}K\| \leq \mu/2$, one has
\begin{align}
\omega\frac{d\phi(\bar{t})}{d\bar{t}}&\leq\|H^{\ast}K\|\sqrt{\frac{p}{q}}+2\|H^{\ast}K\|\phi(\bar{t})+\|H^{\ast}K\|\sqrt{\frac{q}{p}}\phi^{2}(\bar{t})\,. \label{eq:dotPhi_3_dynamic_pf4}
\end{align}
By using the transformation $t =\frac{\bar{t}}{\omega} $, inequality (\ref{eq:dotPhi_3_dynamic_pf4}) and invoking the Comparison Lemma \cite{K:2002}, a lower bound of  the inter-execution time is found as
\begin{align}
\tau^{\ast}=\int_{0}^{1}\dfrac{1}{b_{0}+b_{1}\xi+b_{2}\xi^2+b_{3}\xi^3}{\it d}\xi\,,\label{eq:tauAst_dynamic_pf4}
\end{align}
with $b_{0}=\dfrac{\beta\|H^{\ast}K\|}{\sigma\alpha}$, $b_{1}=2\|H^{\ast}K\|$, $b_{2}=\dfrac{\sigma\alpha\|H^{\ast}K\|}{\beta}$ and $b_{3}=0$.

If $\|H^{\ast}K\| > \mu/2$ and $\gamma \leq 1/(2|H^{\ast}K\|-\mu)$, from (\ref{eq:dotPhi_2_dynamic_pf4}) we get
\begin{align}
\omega &\frac{d\phi(\bar{t})}{d\bar{t}}\leq\|H^{\ast}K\|\sqrt{\frac{p}{q}}+\left(\frac{\mu}{2}+\|H^{\ast}K\|\right)\phi(\bar{t})\nonumber\\
&\quad+\|H^{\ast}K\|\sqrt{\frac{q}{p}}\phi^{2}(\bar{t})+\left(\|H^{\ast}K\|-\frac{\mu}{2}\right)\phi^{3}(\bar{t}) \label{eq:dotPhi_4_dynamic_pf4}
\end{align}
and the minimum dwell-time  $\tau^{\ast}$ satisfies (\ref{eq:tauAst_dynamic_pf4}) with  $b_{0}=\dfrac{\beta\|H^{\ast}K\|}{\sigma\alpha}$, $b_{1}=\dfrac{\mu}{2}+\|H^{\ast}K\|$, $b_{2}=\dfrac{\sigma\alpha\|H^{\ast}K\|}{\beta}$ and $b_{3}=\|H^{\ast}K\|-\dfrac{\mu}{2}$.

Finally, if $\|H^{\ast}K\| > \mu/2$ and $\gamma> 1/(2|H^{\ast}K\|-\mu)$, we obtain
\begin{align}
\omega\frac{d\phi(\bar{t})}{d\bar{t}}&\leq\|H^{\ast}K\|\sqrt{\frac{p}{q}}+\left(2\|H^{\ast}K\|-\frac{1}{2\gamma}\right)\phi(\bar{t})+\nonumber\\
&\quad+\|H^{\ast}K\|\sqrt{\frac{q}{p}}\phi^{2}(\bar{t})+\frac{1}{2\gamma}\phi^{3}(\bar{t})\,, \label{eq:dotPhi_5_dynamic_pf4}
\end{align}
and the lower bound $\tau^{\ast}$ satisfies equation (\ref{eq:tauAst_dynamic_pf4}) with constants $b_{0}=\dfrac{\beta\|H^{\ast}K\|}{\sigma\alpha}$, $b_{1}=2\|H^{\ast}K\|-\dfrac{1}{2\gamma}$, $b_{2}=\dfrac{\sigma\alpha\|H^{\ast}K\|}{\beta}$ and $b_{3}=\dfrac{1}{2\gamma}$. 
Therefore, Zeno behavior is avoided \cite{G:2014}. \hfill $\square$ 

\textcolor{black}{
\section{Periodic Event-Triggered Extremum Seeking} \label{PETC}
}

\textcolor{black}{
Unlike the classical event-triggered extremum seeking considered so far, in periodic event-triggered extremum seeking, the event-triggering condition is verified periodically and at every sampling time it is decided whether or not to compute and to transmit new control signals.}

\textcolor{black}{
This methodology results in two clear advantages, as discussed in \cite{s5}: (a) as the event-triggering condition has to be verified only at the periodic sampling times, instead of continuously, it makes the implementation suitable in standard time-sliced embedded systems architectures; and (b) additionally, the strategy has an inherently guaranteed minimum inter-event time of (at least) one sampling interval of the event-triggering condition, which is easy to tune directly and precludes the possibility of the Zeno behavior naturally.}  

\textcolor{black}{
In a conventional sampled-data feedback setting, the controller would be 
\begin{align}
u(t)=K\hat{G}(t_k) \,, \quad \forall t \in (t_k\,, t_{k+1}]\,,\label{eq:U_sampled_data}
\end{align}
where $t_k$, $k \in \mathbb{N}$, are the sampling times, which are periodic in the sense that $t_k=kh$, for some properly chosen sampling interval $h>0$.}

\textcolor{black}{
Instead of using conventional periodic sampled-data control for extremum seeking, we consider here to use periodic event-triggered control meaning that at each sampling time $t_k=kh$, $k \in \mathbb{N}$, the control values are updated only when certain event-triggering conditions are satisfied. This modifies the controller form (\ref{eq:U_sampled_data}) to  
\begin{align}
u(t)=K\hat{\mathcal{G}}(t) \,, \quad \forall t \in \mathbb{R}_+\,,\label{eq:U_PETC}
\end{align}
where $\hat{\mathcal{G}}(t)$ is a left-continuous signal, given for $t \in (t_k\,, t_{k+1}]$, $k \in \mathbb{N}$, by 
\begin{align}
\hat{\mathcal{G}}(t)=
\begin{cases}
\hat{G}(t_k) \,, \quad \text{when} \quad \mathcal{C}(\hat{G}(t_k) \,,\hat{\mathcal{G}}(t_k))>0\,,\\
\hat{\mathcal{G}}(t_k) \,, \quad \text{when} \quad \mathcal{C}(\hat{G}(t_k) \,,\hat{\mathcal{G}}(t_k))\leq 0\,,
\end{cases}
\label{eq:left_continuous}
\end{align}
and some initial value for $\hat{\mathcal{G}}(0)$. Hence, the value of $\hat{\mathcal{G}}(t)$ can be interpreted as the most recently transmitted measurement of the gradient estimate $\hat{G}$ to the controller at time $t$. Whether or not new gradient measurements are transmitted to the controller is based on the event-triggering condition $\mathcal{C}$. In particular, if at time $t_k$ it holds that $\mathcal{C}(\hat{G}(t_k) \,,\hat{\mathcal{G}}(t_k))>0$, the signal $\hat{G}(t_k)$ is transmitted through the controller and $\hat{\mathcal{G}}$ as well as the control value $u$ are updated accordingly. In case $\mathcal{C}(\hat{G}(t_k) \,,\hat{\mathcal{G}}(t_k))\leq 0$, no new gradient information is sent to the controller, in which case the input $u$ is is not updated and kept the same for (at least) another sampling interval implying that no control computations are needed and no new gradient measurements and control values have to be transmitted. }

\textcolor{black}{
In this section, we restrict ourselves to the same event-triggering conditions (static or dynamic) used earlier for classical event-triggered extremum seeking, by choosing $\mathcal{C}(\hat{G}(t_k) \,,\hat{\mathcal{G}}(t_k)):=\Xi(\hat{G},e)$, with 
\begin{align}
\Xi(\hat{G},e)&=\sigma \alpha \|\hat{G}(t)\|^2-\beta \|e(t)\|\|\hat{G}(t)\|	\,, \label{eq:Xi_event_1_PETC}
\end{align}
given in (\ref{eq:Xi_event_1}), but redefining the error vector as 
\begin{align}
e(t):=\hat{\mathcal{G}}(t_{k})-\hat{G}(t) \,, \quad \forall t \in \lbrack t_{k}\,, t_{k+1}) \,, \quad k\in \mathbb{N} \,. \label{eq:e_event_PETC}
\end{align}
The analysis can be carried out using the step-by-step presented in \cite{s5}.}

\section{Simulation Results} \label{sec:sim}

We consider  the multivariable nonlinear map (\ref{eq:Q_1_event}) with an input $\theta(t)\in \mathbb{R}^{2}$, an output $y(t) \in \mathbb{R}$, and unknown parameters
\begin{align}
H=\begin{bmatrix} 100 & 30 \\ 30 & 20\end{bmatrix}>0\,,
\end{align} 
$Q^{\ast}=100$ and $\theta^{\ast}=\begin{bmatrix}2 & 4\end{bmatrix}^{T}$. 
The dither vectors (\ref{eq:S_event}) and (\ref{eq:M_event}) have parameters $a_{1}=a_{2}=0.1$, $\omega_{1}=0.1$ [rad/sec], and $\omega_{2}=0.7$ [rad/sec], as in \cite{GKN:2012}, and we select the event-triggered parameters $\sigma=0.5$, $\alpha=1$, $\beta=3.1521$, $\mu=0.4320$ and $\gamma=0.0542$. The control gain matrix is $K=10^{-2}\begin{bmatrix} -6 & 0 \\ 0 & -20\end{bmatrix}$ and initial condition is $\upsilon(0)=0$. \textcolor{black}{Due to space limitations, we are going to present numerical simulations only to the case of classical event-triggered extremum seeking designs.}

\begin{figure}[h!]
	\centering
	\subfigure[\underline{\bf Static:} sample-and-hold gradient estimate, $\hat{G}(t_{k})$.\label{fig:S_sigma_0_5_hatG_tk}]{\includegraphics[width=4.1cm]{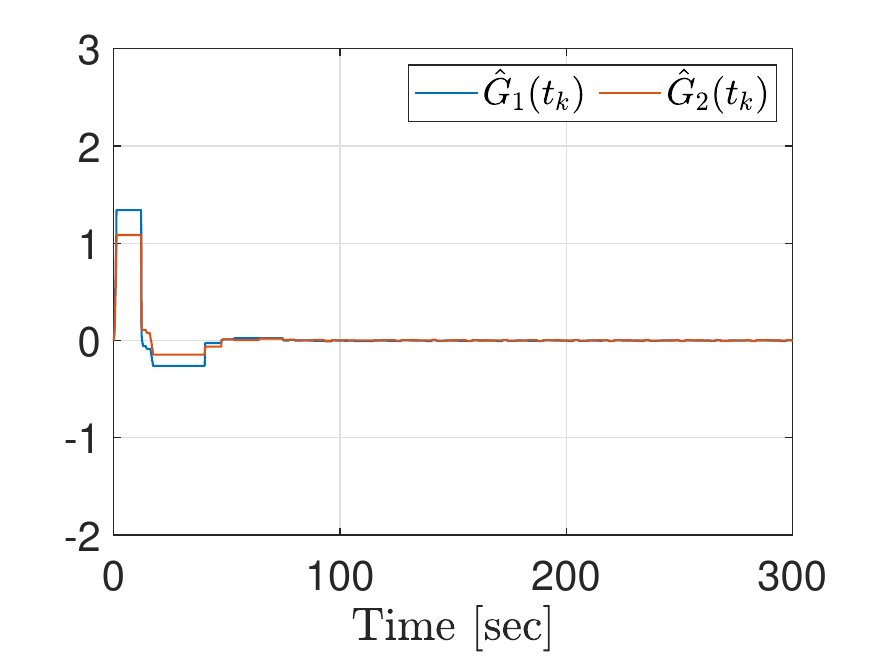}}
	~
	\subfigure[\underline{\bf Dynamic:} sample-and-hold gradient estimate, $\hat{G}(t_{k})$. \label{fig:D_sigma_0_5_hatG_tk}]{\includegraphics[width=4.1cm]{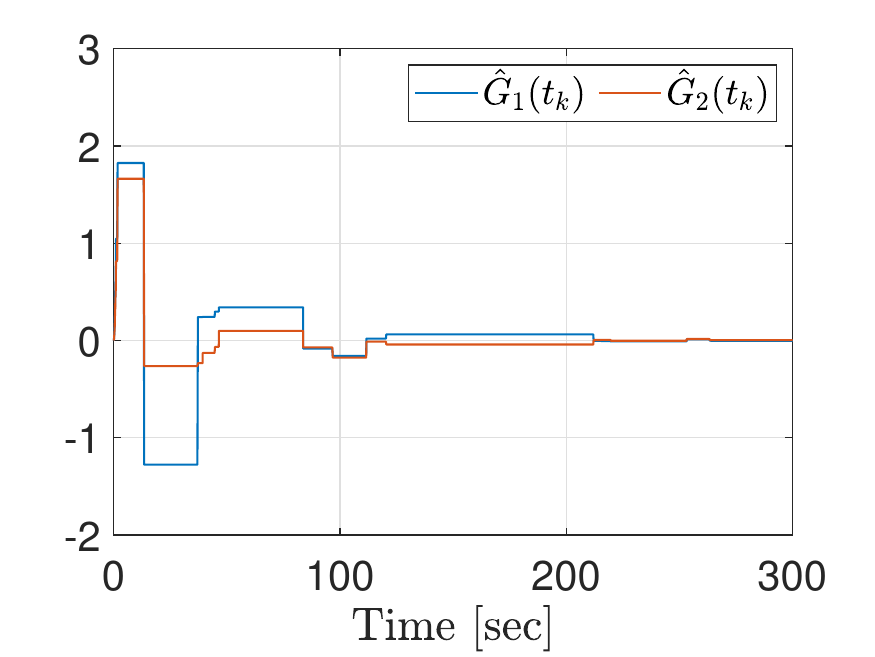}}
	\\
	\subfigure[\underline{\bf Static:} control input, $U(t)$. \label{fig:S_sigma_0_5_U}]{\includegraphics[width=4.1cm]{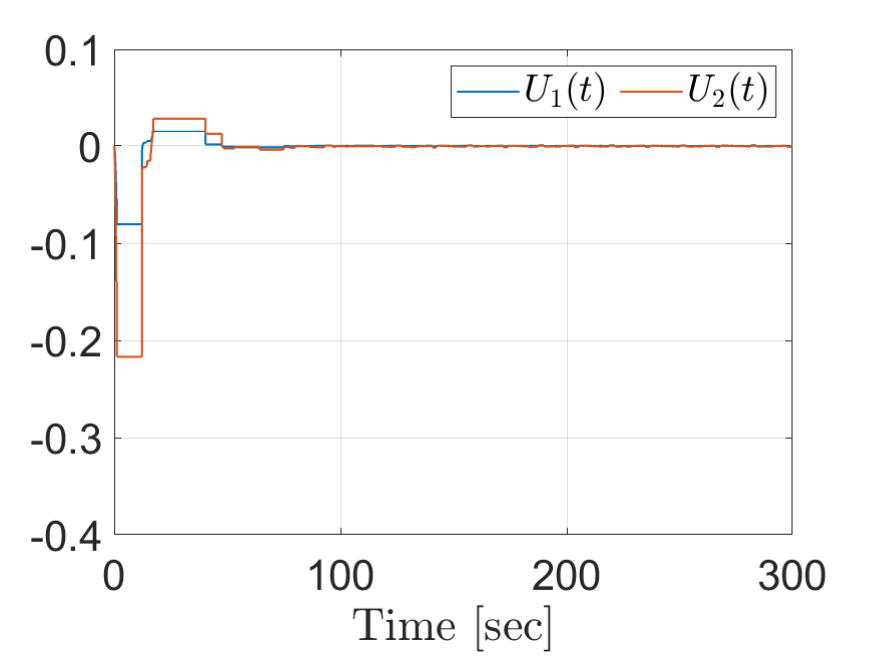}}
	~
	\subfigure[\underline{\bf Dynamic:} control input, $U(t)$. \label{fig:D_sigma_0_5_U}]{\includegraphics[width=4.1cm]{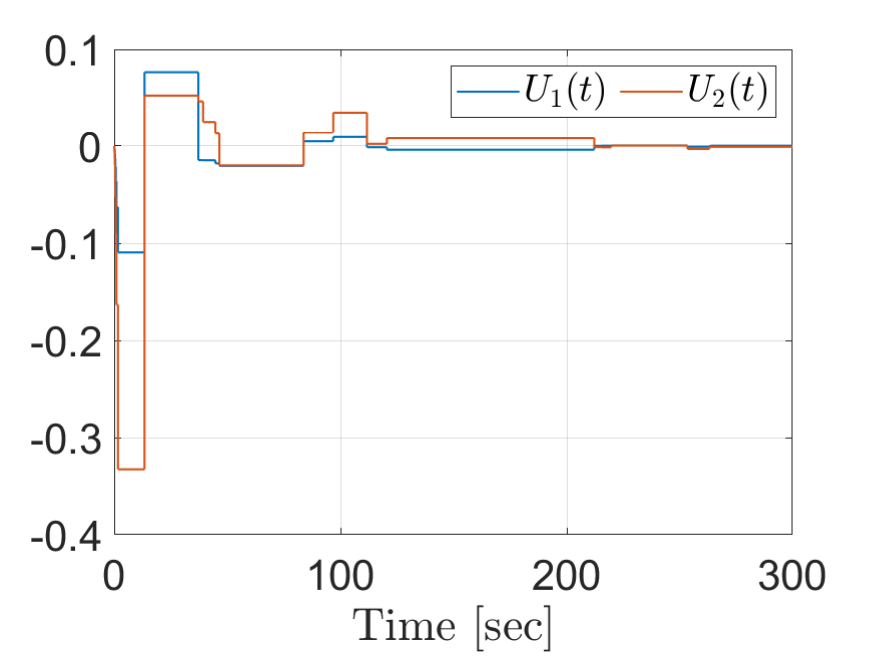}}
	\\
	\subfigure[\underline{\bf Static:} input of the nonlinear map, $\theta(t)$. \label{fig:S_sigma_0_5_theta}]{\includegraphics[width=4.1cm]{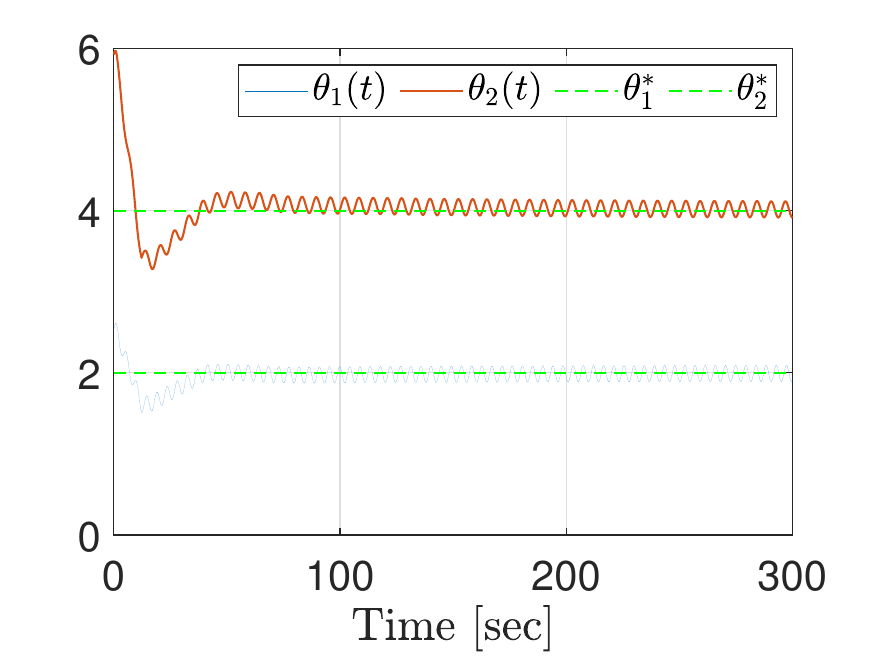}}
	~
	\subfigure[\underline{\bf Dynamic:} input of the nonlinear map, $\theta(t)$. \label{fig:D_sigma_0_5_theta}]{\includegraphics[width=4.1cm]{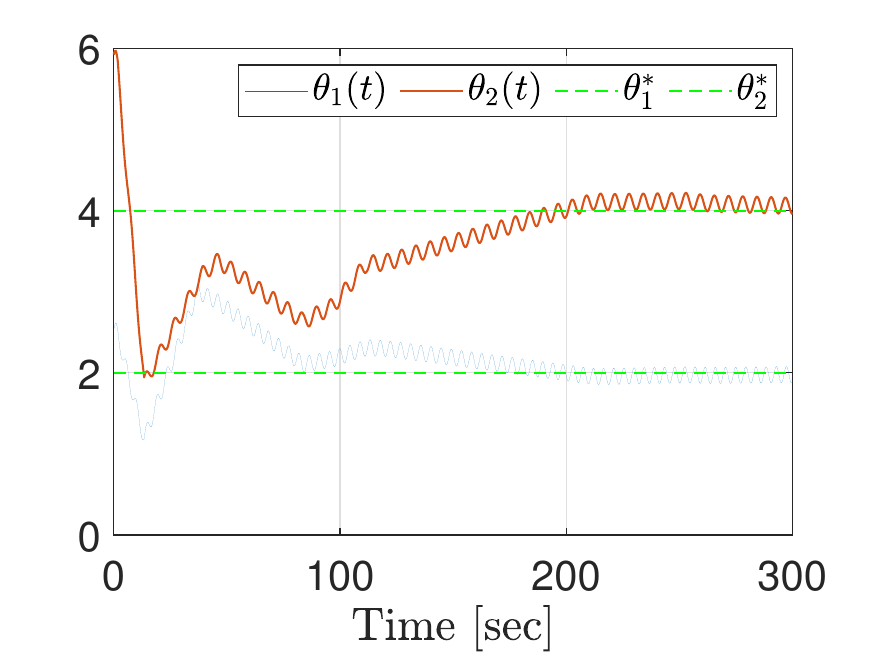}}
	\\
	\subfigure[\underline{\bf Static:} output of the nonlinear map, $y(t)$. \label{fig:S_sigma_0_5_y}]{\includegraphics[width=4.1cm]{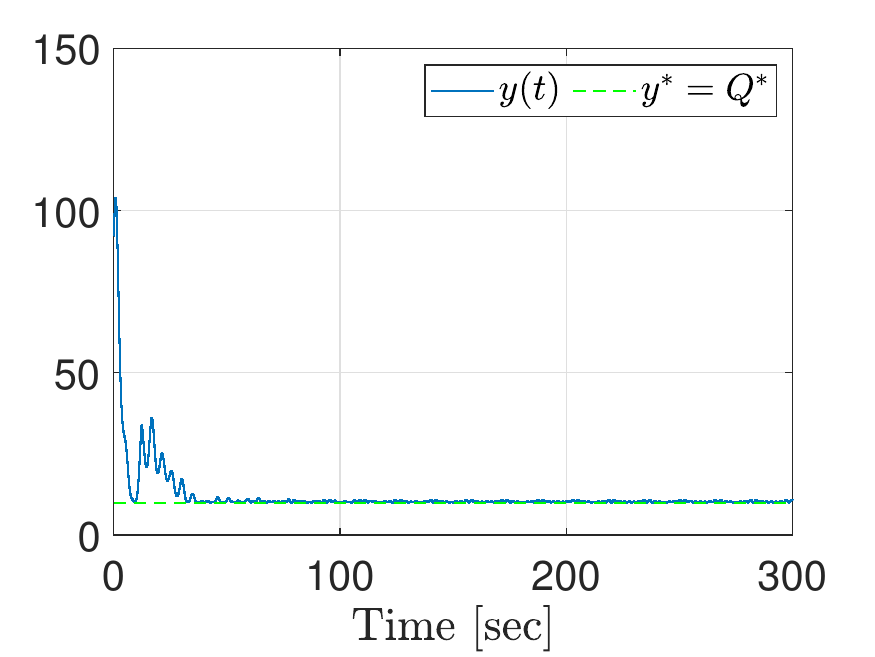}}
	~
	\subfigure[\underline{\bf Dynamic:} output of the nonlinear map, $y(t)$. \label{fig:D_sigma_0_5_y}]{\includegraphics[width=4.1cm]{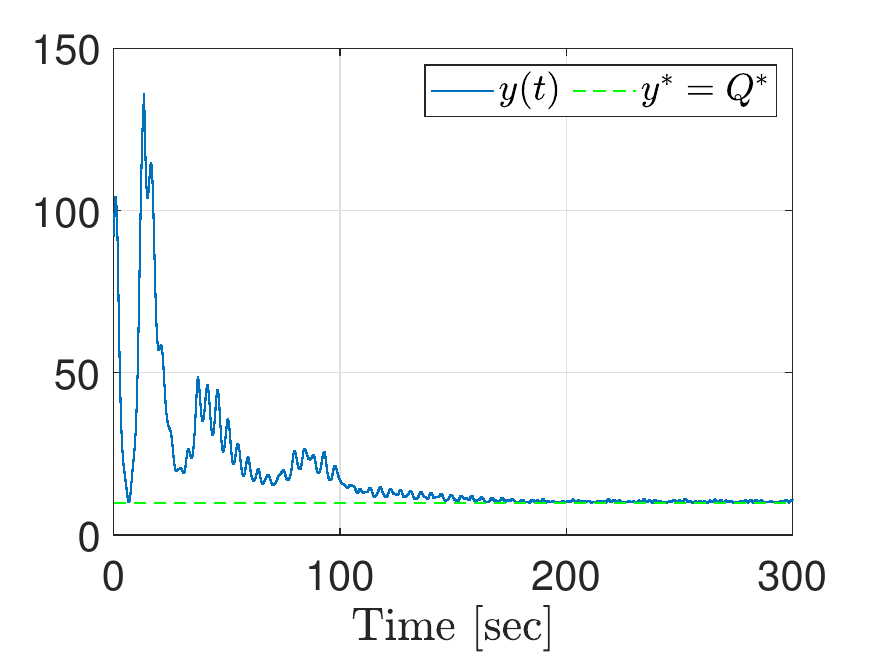}}
	\caption{Static and Dynamic Event-triggered Extremum Seeking Systems. \label{fig:SD_ETESC}}
\end{figure}

In Fig.~\ref{fig:SD_ETESC}, both static and dynamic event-triggered approaches are simulated, with initial condition $\hat{\theta}(0)=\begin{bmatrix}2.5, & 6\end{bmatrix}^{T}$. Fig.~\ref{fig:S_sigma_0_5_hatG_tk} and Fig.~\ref{fig:D_sigma_0_5_hatG_tk} shows the convergence to zero  of both the sampled-and-hold version of the gradient estimate   when the control signals are given by Fig.~\ref{fig:S_sigma_0_5_U} and Fig.~\ref{fig:D_sigma_0_5_U}, respectively. Of course, the gradient stabilization implies  reaching the optimizer $\theta^{\ast}$, as illustrated in Fig.~\ref{fig:S_sigma_0_5_theta} and Fig.~\ref{fig:D_sigma_0_5_theta}, consequently, the variable $y(t)$ reaches its extremum value as shown in Fig.~\ref{fig:S_sigma_0_5_y} and Fig.~\ref{fig:D_sigma_0_5_y}. 

\begin{figure}[h!]
	\centering
  \subfigure[Component of the control signal, $U_{1}(t)$. \label{fig:3D_D_sigma_0_5_U1}]{\includegraphics[width=4.1cm]{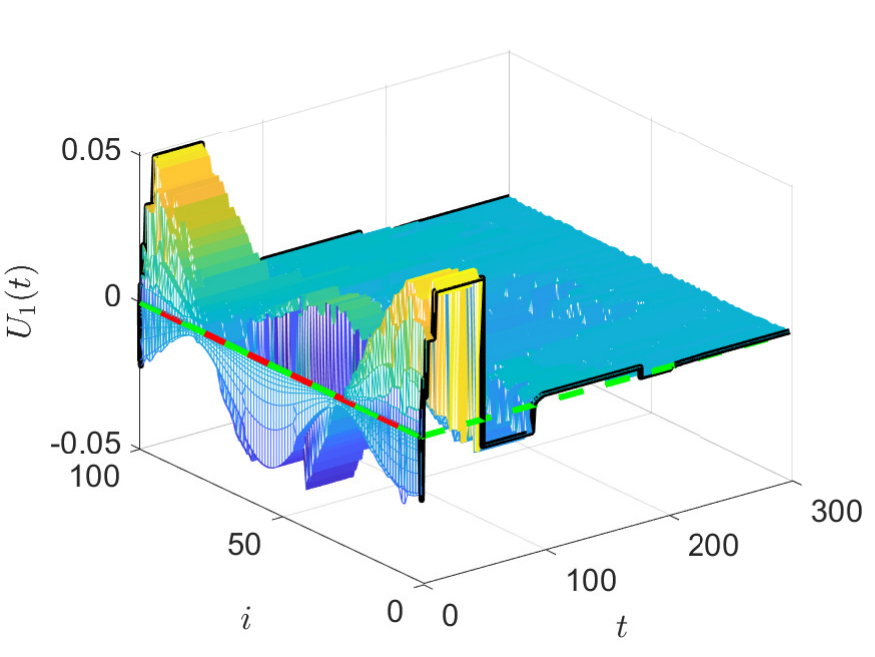}}
	~
	\subfigure[Component of the control signal, $U_{2}(t)$. \label{fig:3D_D_sigma_0_5_U2}]{\includegraphics[width=4.1cm]{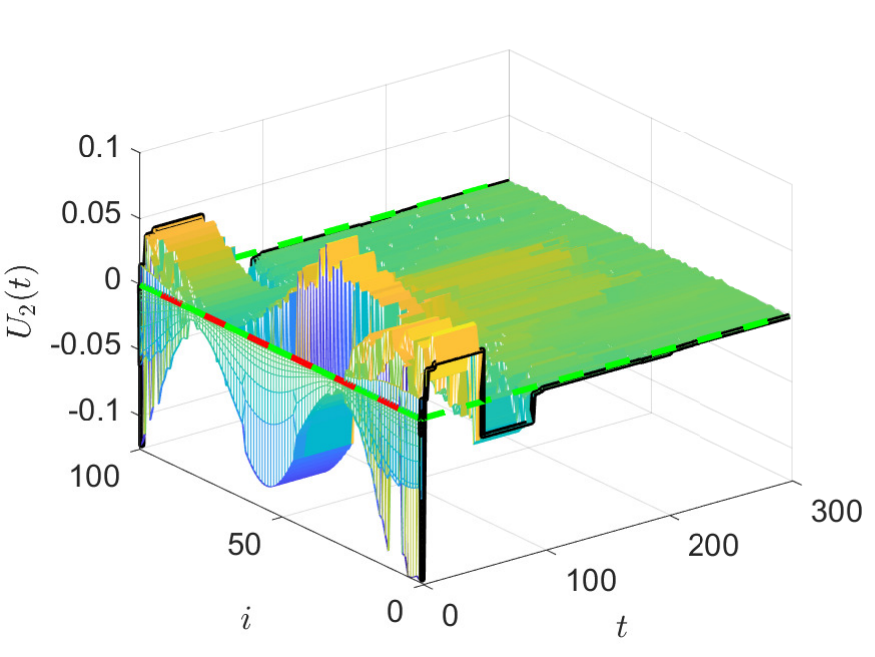}}
	\\	
	\subfigure[Component of the gradient estimate, $\hat{G}_{1}(t)$. \label{fig:3D_D_sigma_0_5_hatG1}]{\includegraphics[width=4.1cm]{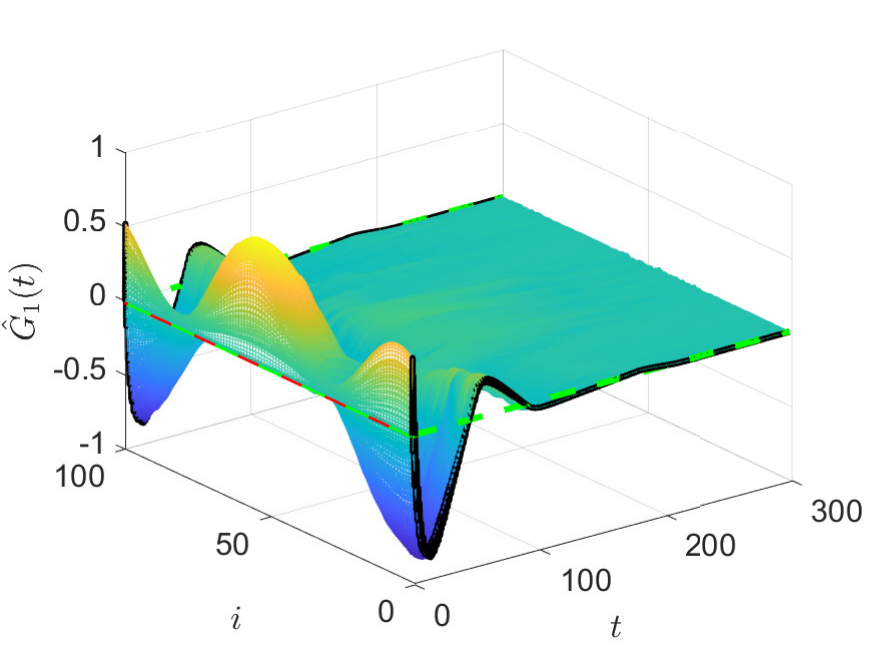}}
	~
	\subfigure[Component of the gradient estimate, $\hat{G}_{2}(t)$. \label{fig:3D_D_sigma_0_5_hatG2}]{\includegraphics[width=4.1cm]{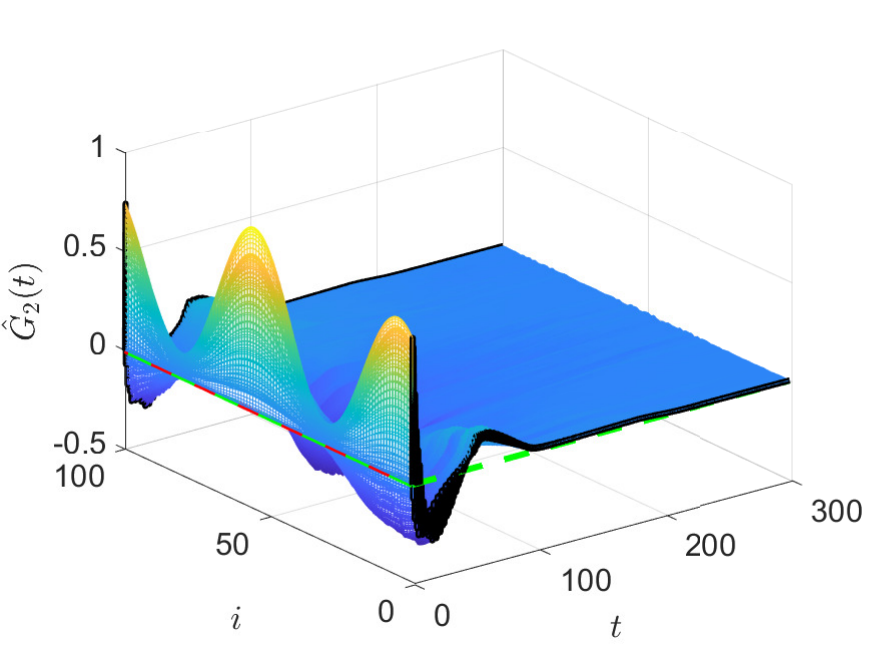}}
	\\
	\subfigure[Component of the nonlinear map input, $\theta_{1}(t)$. \label{fig:3D_D_sigma_0_5_theta1}]{\includegraphics[width=4.1cm]{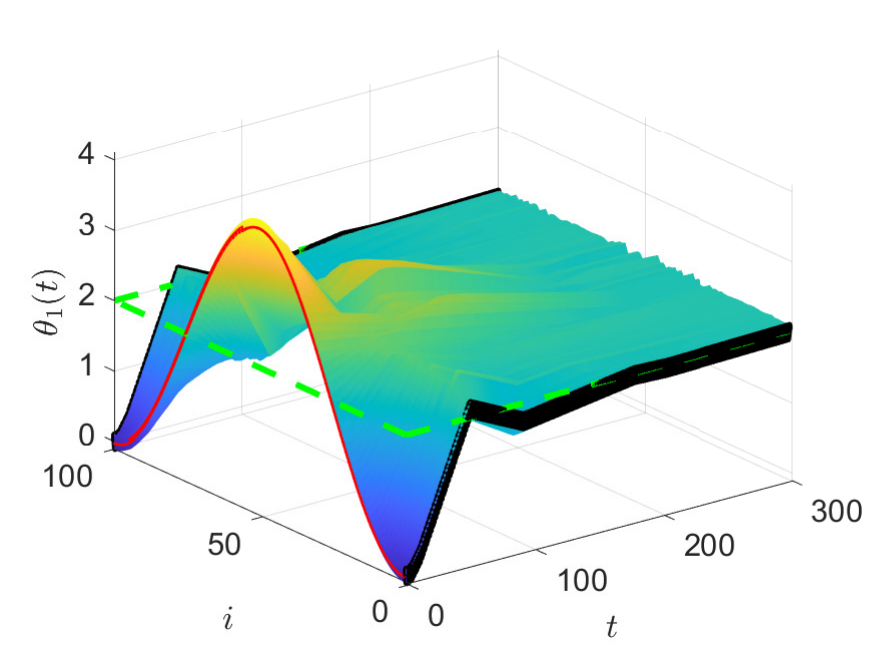}}
	~
	\subfigure[Component of the nonlinear map input, $\theta_{2}(t)$. \label{fig:3D_D_sigma_0_5_theta2}]{\includegraphics[width=4.1cm]{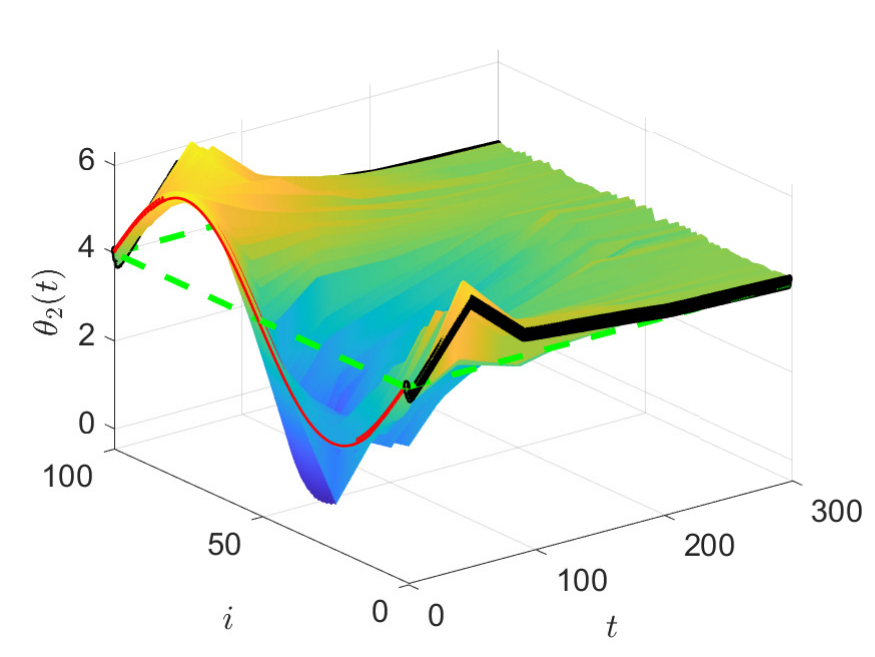}}
	\caption{Dynamic Event-triggered Extremum Seeking Feedback System. \label{fig:3D_D_sigma_0_5}}
\end{figure}

\begin{figure}[h!]
\centering
\includegraphics[width=7.5cm]{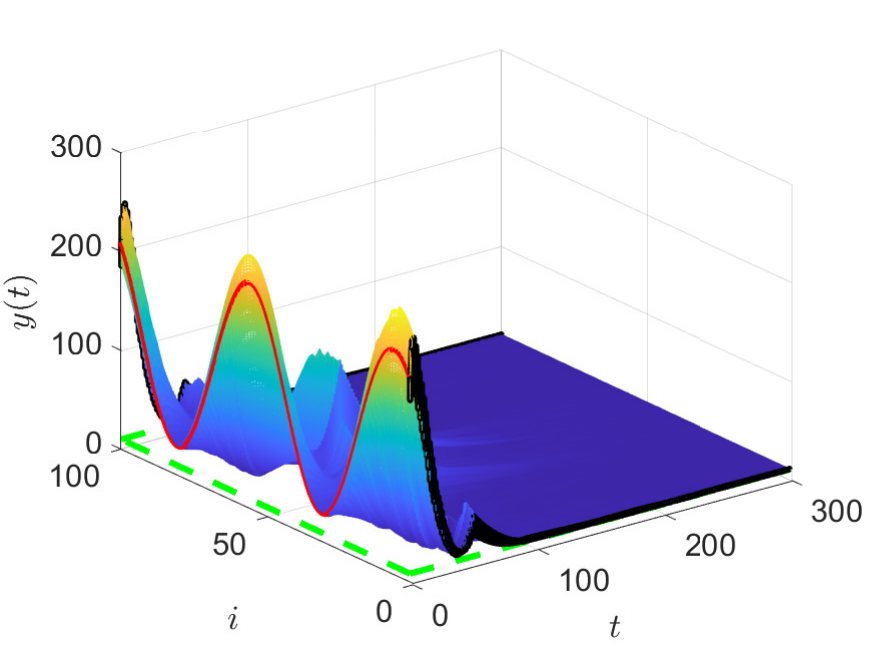}
\caption{Output of the nonlinear map, $y(t)$. \label{fig:3D_D_sigma_0_5_y}} 
\end{figure}

Hereafter, we have considered several simulations by using the set of initial conditions $$\hat{\theta}(0)=\begin{bmatrix}2-2\cos \left(\frac{2\pi}{100}i\right), & 4-2\sin \left(\frac{2\pi}{100}i\right)\end{bmatrix}^{T}$$ for $i=1,\ldots,100$. Fig.~\ref{fig:3D_D_sigma_0_5} shows the time-evolution of the proposed dynamic event-triggered extremum seeking approaches. For all initial conditions, the input signal in Fig.~\ref{fig:3D_D_sigma_0_5_U1} and Fig.~\ref{fig:3D_D_sigma_0_5_U2} ensures  convergence of the gradient estimate to neighborhood of zero 
as presented in  
Fig.~\ref{fig:3D_D_sigma_0_5_hatG1} and Fig.~\ref{fig:3D_D_sigma_0_5_hatG2}. Thus, $\theta_1(t)$ and $\theta_{2}(t)$ tend asymptotically to $\theta_1^{\ast}$ and $\theta_{2}^{\ast}$, respectively, consequently, $y(t)$ to $Q^{\ast}$ \textcolor{black}{(see Fig.~\ref{fig:3D_D_sigma_0_5_theta1},  Fig.~\ref{fig:3D_D_sigma_0_5_theta2}, and  Fig.~\ref{fig:3D_D_sigma_0_5_y}).} 

\vspace{-0.5cm}

Finally, Table~\ref{tab:IEI} summarizes the statistical data obtained for a set of 5400 simulations of 300 seconds. For any value $\sigma$, it is possible to verify that the interval between executions as well as the dispersion measures (mean deviation, variance and standard deviation) are greater when the dynamic strategy is employed. For instance, on average, with $\sigma=0.001$, the static strategy performs 1021 updates while the dynamic one needs 755 updates. When, $\sigma = 0.9$, the static case requires 46 updates versus 15 for the dynamic case. Although the control objective is achieved with both strategies, it is worth noting that the dynamic strategy achieves the extremum employing less control effort and requiring a small number of control updates when compared with the static approach (see Table~\ref{tab:IEI}).

\begin{table*}[h!]
\centering
\caption{Statistics of the inter-execution intervals, $t_{k+1}-t_{k}$. \label{tab:IEI}}
\begin{scriptsize}
\begin{tabular}{|l|c|c|c|c|c|c|c|c|c|c|}
\cline{2-9} 
 \multicolumn{1}{c}{}& \multicolumn{2}{|c|}{Mean} & \multicolumn{2}{|c|}{Mean Deviation} & \multicolumn{2}{|c|}{Variance} & \multicolumn{2}{|c|}{Standard Deviation} \\
\hline
\multicolumn{1}{|c|}{$\sigma$} & Static & Dynamic & Static & Dynamic & Static & Dynamic & Static & Dynamic \\
\hline
0.001    & 0.2939 &  0.3974 &  0.5088 &  0.7330 &  13.9295 &  22.3420 &  3.7322 &  4.7268 \\
0.002    & 0.3099 &  0.4381 &  0.5337 &  0.8048 &  14.3554 &  23.8753 &  3.7888 &  4.8862 \\
0.003    & 0.3301 &  0.7370 &  0.5670 &  1.3255 &  15.5122 &  39.9080 &  3.9386 &  6.3173 \\
0.004    & 0.3496 &  0.7096 &  0.5983 &  1.2811 &  16.2814 &  38.8355 &  4.0350 &  6.2318 \\
0.005    & 0.3818 &  0.8633 &  0.6512 &  1.5367 &  17.6344 &  48.3978 &  4.1993 &  6.9569 \\                
0.006    & 0.4150 &  0.8587 &  0.7029 &  1.5310 &  18.5522 &  48.7397 &  4.3072 &  6.9814 \\
0.007    & 0.4365 &  1.0187 &  0.7369 &  1.7873 &  19.3296 &  56.6231 &  4.3965 &  7.5248 \\
0.008    & 0.4260 &  0.7903 &  0.7179 &  1.4177 &  18.8293 &  43.4613 &  4.3393 &  6.5925 \\
0.009    & 0.4284 &  0.8496 &  0.7209 &  1.5189 &  19.0454 &  48.5479 &  4.3641 &  6.9676 \\
0.010    & 0.4686 &  0.9245 &  0.7843 &  1.6425 &  20.8145 &  53.0984 &  4.5623 &  7.2869 \\
0.020    & 0.6191 &  1.3454 &  1.0229 &  2.3079 &  28.7969 &  73.6432 &  5.3663 &  8.5816 \\
0.030    & 0.7958 &  1.4215 &  1.2801 &  2.4192 &  37.2480 &  78.3938 &  6.1031 &  8.8540 \\
0.040    & 0.8192 &  1.5215 &  1.3319 &  2.5772 &  40.4228 &  85.4836 &  6.3579 &  9.2457 \\
0.050    & 0.9848 &  1.9907 &  1.5642 &  3.2645 &  48.6906 & 109.4917 &  6.9779 & 10.4638 \\
0.060    & 1.1553 &  2.4570 &  1.8035 &  3.9291 &  55.7665 & 133.7192 &  7.4677 & 11.5637 \\
0.070    & 1.2836 &  2.3071 &  1.9779 &  3.7365 &  65.7437 & 127.1205 &  8.1082 & 11.2748 \\
0.080    & 1.1995 &  2.8969 &  1.8986 &  4.5692 &  61.9130 & 160.1954 &  7.8685 & 12.6568 \\
0.090    & 1.4664 &  2.9279 &  2.2550 &  4.5974 &  74.3021 & 163.4882 &  8.6199 & 12.7862 \\ 
0.100    & 1.2588 &  2.5863 &  1.9952 &  4.1467 &  64.1229 & 147.5280 &  8.0077 & 12.1461 \\
0.200    & 1.9571 &  5.1650 &  3.1000 &  7.7201 & 102.7835 & 286.3306 & 10.1382 & 16.9213 \\  
0.300    & 2.0788 &  7.0843 &  3.3859 & 10.2374 & 113.2025 & 391.6320 & 10.6397 & 19.7897 \\ 
0.400    & 3.4475 &  9.3378 &  5.3195 & 12.8145 & 180.9892 & 479.3062 & 13.4532 & 21.8931 \\
0.500    & 3.9839 & 10.5757 &  6.1247 & 14.4068 & 218.1339 & 566.1838 & 14.7694 & 23.7946 \\
0.600    & 3.4004 & 12.8756 &  5.5095 & 16.9665 & 191.6231 & 608.7021 & 13.8428 & 24.6719 \\
0.700    & 5.0007 & 14.7065 &  7.6219 & 18.9248 & 259.6896 & 677.2911 & 16.1149 & 26.0248 \\
0.800    & 5.6172 & 16.5772 &  8.4851 & 21.8088 & 297.7478 & 807.1438 & 17.2554 & 28.4103 \\
0.900    & 6.6610 & 20.2129 & 10.1423 & 25.3364 & 445.3767 & 876.5687 & 21.1039 & 29.6069 \\
\hline
\end{tabular}
\end{scriptsize}
\end{table*}

\section{Conclusion} \label{sec:concl}

This paper brings the first effort to pursue a combination of event-triggered control and extremum seeking. We proposed the static and dynamic multivariable event-triggered extremum seeking based on the perturbation method for nonlinear static maps. The contribution of treating the specific hybrid learning dynamics \cite{PT:2017}  constituted by the sample-and-hold mechanism and the asynchronous nature of the resulting ``sampled-data system'' for event-triggering extremum seeking actions based on averaging theory and periodic perturbation method is clear. The approach provides explicit functional forms for the exponential convergence characterizing \textcolor{black}{the stability properties and Zeno preclusion for the closed-loop average system. This allows us to effectively optimize nonlinear maps in real-time without relying on prior knowledge of their parameters, with guaranteed practical  exponential stability to the non-average (real) system and convergence to a small neighborhood of the extremum point.}   
%
%
\textcolor{black}{Alternatively, the Zeno behavior can also be avoided when we consider periodic event-triggered architecture \cite{s5}, since accumulations of update times cannot occur by construction.}

\textcolor{black}{Although we kept the presentation using the classical Gradient extremum seeking, the generalization of the results for the Newton-based approach is also possible \cite{GKN:2012}. The extension of the proposed approach for time-varying maps \cite{PKB:2023}, where unexpected events occurring at unexpected times which modify the extremum point seems to be a challenging problem to be considered as well.} Finally, the proposed approach can also be expanded to PDE systems following the ideas introduced in  \cite{FKKO:2018,REOGK:2021,REOGK:2020,OK:2021,OK:2022_book_golden,OKT:2017} and 
 \cite{RDEK:2021,rathnayake2022sampled}.

\begin{small}

\section*{Appendix} 
\appendix

\textcolor{black}{It is important to clarify that
the following averaging result for differential inclusions with discontinuous right-hand sides by Plotnikov \cite{P:1979} takes into account discontinuities not in the periodic perturbations, but in the states. In other words, in our case, we still have periodic perturbations even if the discontinuities are not periodic in the triggering events. Hence, we are indeed able to invoke it in our analysis carried out in the proofs of Theorems \ref{thm:NETESC_2} and \ref{thm:NETESC_4}.}

\section{Average of Discontinuous Systems}
\label{appendix_plotnikov}

From \cite{P:1979}, let us consider the differential inclusion
\begin{align}
\frac{d x}{dt} \in \varepsilon X(t,x)\,, \quad x(0)=x_{0}\,, \label{eq:A1}
\end{align}
where $x$ is an n-dimensional vector, $t$ is time, $\varepsilon$ is a small parameter, and $X(t,x)$ is a multivalued function that is $T$-periodic in $t$ and puts in correspondence with with each point $(t,x)$ of a certain domain of the ($n+1$)-dimensional space a compact set $X(t,x)$ of the $n$-dimensional space.

Let us put in correspondence with the inlusion (\ref{eq:A1}) the average inclusion
\begin{align}
\frac{d \xi}{dt} \in \varepsilon \bar{X}(\xi)\,, \quad \xi(0)=x_{0}\,, \label{eq:A2}
\end{align}
where
\begin{align}
\bar{X}(\xi)=\frac{1}{T}\int_{0}^{T}X(\tau,\xi)d\tau\,. \label{eq:A3}
\end{align}

\begin{theorem} \label{thm:A1}
Let a multivalued mapping $X(t,x)$ be defined in the domain $Q\left\{t\geq 0\,, x\in D\subset\mathbb{R}^{n}\right\}$ and let in this domain the set $X(t,x)$ be a nonempty compactum for all admissible values of the arguments and the mapping $X(t,x)$ be continuous and uniformly bounded and satisfy the Lipschitz condition with respect to $x$ with a constant $\lambda$, {\it i.e.}, $X(t,x) \subset S_{M}(0)$, $\delta(X(t,x')-X(t,x''))\leq \lambda \|x'-x''\|$, where $\delta(P,Q)$ is the Hausdorff distance between the sets $P$ and $Q$, {\it i.e.}, $\delta(P,Q)=\min\left\{d|P\subset S_{d}(Q), Q \subset S_{d}(P)\right\}$, $S_{d}(N)$ being the d-neighborhood of a set $N$ in the space $\mathbb{R}^{n}$; the mapping $X(t,x)$ be $T$-periodic in $t$; for all $x_{0}\in D'\subset D$ the solutions of inclusion (\ref{eq:A2}) lie in  the domain $D$ together with a certain $\rho$-neighborhood. Them for each $L>0$ there exist $\varepsilon^{0}(L)>0$ and $c(L)>0$ such that for $\varepsilon \in \rbrack0\,,\varepsilon^{0}\rbrack$ and $t \in \lbrack 0, L\varepsilon^{-1}\rbrack$:
\begin{enumerate}
\item for each solution $x(t)$ of the inclusion (\ref{eq:A1}) there exists a solution $\xi(t)$ of the inclusion (\ref{eq:A2}) such that 
\begin{align}
\|x(t)-\xi(t)\|\leq c\varepsilon =\mathcal{O}(\varepsilon); \label{eq:A4}
\end{align}
\item for each solution $\xi(t)$ of the inclusion (\ref{eq:A2}) there exists a solution $x(t)$ of the inclusion (\ref{eq:A1}) such that  the inequality (\ref{eq:A4}) holds.
\end{enumerate}
Thus the following estimate is valid:
\begin{align}
\delta(\bar{R}(t),R'(t))\leq c\varepsilon =\mathcal{O}(\varepsilon)\,, \label{eq:A5}
\end{align}
where $\bar{R}(t)$ is a section of the family of solutions of the inclusion (\ref{eq:A2}) and $R'(t)$ is the closure of the section $R(t)$ of the family of solutions of the inclusion (\ref{eq:A1}).
\end{theorem}

\textit{Proof:} The proof of Theorem~\ref{thm:A1} is constructed following Banfy's theorem \cite{B:1967} and  replacing references in the first Bogolyubov's theorem with the ones in \cite{P:1979}. 
\hfill $\square$

\begin{theorem} \label{thm:A2}
Let all the conditions of Theorem~\ref{thm:A1} and also the following condition be fulfilled: the $R$-solution $\bar{R}(t)$ of inclusion (\ref{eq:A2}) is uniformly asymptotically stable. Then there exist $\varepsilon^{0}>0$ and $c>0$ such that for $0< \varepsilon \leq \varepsilon^{0}$
\begin{align}
\delta(\bar{R}(t),R'(t))\leq c\varepsilon = \mathcal{O}(\varepsilon)\,, \label{eq:A6}
\end{align}
for all $t\geq 0$.
\end{theorem}

\textit{Proof:} The proof of Theorem~\ref{thm:A2} can be derived by straightforwardly following  \cite{F:1971}. Substituting the references to Bogolyubov's theorem to Theorem~\ref{thm:A1}, all estimates in the sense of the Hausdorff metric can be found in reference \cite{P:1979}. \hfill $\square$ 

\end{small}


\renewcommand{\baselinestretch}{0.85}

\end{document}